\pdfoutput=1
\documentclass[review]{elsarticle}
\usepackage{graphicx,epsf,subfigure,pstricks,pst-node,psfrag}
\usepackage[colorlinks]{hyperref}
\hypersetup{
 citecolor={red}
}
\usepackage{soul}
\usepackage{amsthm,amssymb,amsmath}
\usepackage{algorithm}
\usepackage{algpseudocode}
\usepackage{booktabs}
\usepackage{cleveref}
\crefname{equation}{}{}
\usepackage{multirow}

\algblock{Input}{EndInput}
\algnotext{EndInput}
\algblock{Output}{EndOutput}
\algnotext{EndOutput}

\newcommand{\myRe}{\mathbb{R}}
\newcommand{\myEx}{\mathbb{E}}

\newcommand{\bal}{\boldsymbol{\alpha}}
\newcommand{\bbe}{\boldsymbol{\beta}}

\crefname{figure}{Figure}{Figures}
\Crefname{figure}{Figure}{Figures}
\crefname{equation}{}{}
\Crefname{equation}{}{}
\crefname{appsec}{}{}
\usepackage{bm}

  \def\clap#1{\hbox to 0pt{\hss#1\hss}}

\providecommand{\mat}[1]{\bm{#1}}%
\renewcommand{\vec}[1]{\mathbf{#1}}

\providecommand{\mA}{\ensuremath{\mat{A}}}
\providecommand{\mB}{\ensuremath{\mat{B}}}

\providecommand{\mE}{\ensuremath{\mat{E}}}

\providecommand{\mI}{\ensuremath{\mat{I}}}

\providecommand{\mM}{\ensuremath{\mat{M}}}

\providecommand{\mQ}{\ensuremath{\mat{Q}}}
\providecommand{\mR}{\ensuremath{\mat{R}}}

\providecommand{\mU}{\ensuremath{\mat{U}}}

\providecommand{\mW}{\ensuremath{\mat{W}}}

\providecommand{\ve}{\ensuremath{\vec{e}}}
\providecommand{\vf}{\ensuremath{\vec{f}}}

\providecommand{\vi}{\ensuremath{\vec{i}}}
\providecommand{\vj}{\ensuremath{\vec{j}}}

\providecommand{\vm}{\ensuremath{\vec{m}}}

\providecommand{\vp}{\ensuremath{\vec{p}}}
\providecommand{\vq}{\ensuremath{\vec{q}}}
\providecommand{\vr}{\ensuremath{\vec{r}}}
\providecommand{\vs}{\ensuremath{\vec{s}}}

\providecommand{\vu}{\ensuremath{\vec{u}}}

\providecommand{\vw}{\ensuremath{\vec{w}}}
\providecommand{\vx}{\ensuremath{\vec{x}}}
\providecommand{\vy}{\ensuremath{\vec{y}}}

\begin{document}
\title{\textbf{Extremum Sensitivity Analysis with Polynomial Monte Carlo Filtering} }
\date{}

\author{Chun Yui Wong\corref{cor1}%
\fnref{fn1}}
\ead{cyw28@cam.ac.uk}
\author{Pranay Seshadri\fnref{fn2}}

\author{Geoffrey Parks\fnref{fn3}}

\cortext[cor1]{Corresponding author}
\fntext[fn1]{PhD student, Department of Engineering, University of Cambridge.}
\fntext[fn2]{Research Fellow, Department of Mathematics (Statistics Section), Imperial College London.}
\fntext[fn3]{Reader, Department of Engineering, University of Cambridge.}

\begin{abstract}

Global sensitivity analysis is a powerful set of ideas and heuristics for understanding the importance and interplay between uncertain parameters in a computational model. Such a model is characterized by a set of input parameters and an output quantity of interest, where we typically assume that the inputs are independent and their marginal densities are known. If the output quantity is smooth, polynomial chaos can be used to extract Sobol’ indices.

In this paper, we build on these well-known ideas by examining two different aspects of this paradigm. First, we study whether sensitivity indices can be computed efficiently if one leverages a polynomial ridge approximation—a polynomial fit over a subspace. Given the assumption of anisotropy in the dependence of a function, we show that sensitivity indices can be computed with a reduced number of model evaluations. Second, we discuss methods for evaluating sensitivities when constrained near output extrema. Methods based on the analysis of skewness are reviewed and a novel type of indices based on Monte Carlo filtering (MCF)---extremum Sobol’ indices---is proposed. We combine these two ideas by showing that these indices can be computed efficiently with ridge approximations, and explore the relationship between  MCF-based indices and skewness-based indices empirically.
\end{abstract}


\maketitle
\setcounter{footnote}{0}

\section{Introduction and motivation}
Across all data-centric disciplines, sensitivity analysis is an important task. The phrase \emph{``which of my parameters is the most important?”} can be heard across shop floors, office spaces and board rooms. In global sensitivity analysis, the measure of importance is typically the conditional variance, leading to the familiar Sobol' indices \cite{sobol1993sensitivity,sobol2001global}, the related total indices \cite{homma1996importance} and the Fourier amplitude sensitivity test (FAST) indices \cite{cukier1978nonlinear}. However, these are not the only metrics available: over the years, other metrics for sensitivity analysis have been proposed, including those based on conditional skewness \cite{wang2001quantifying} and kurtosis \cite{delloca2017moment-based}. Beyond this, Hooker \cite{hooker2004discovering} and Liu and Owen\cite{liu2006estimating} introduce the notion of a superset importance, which is the sum of all Sobol' indices based on a superset of an input set. Techniques based on partial derivative evaluations of a function have also been extensively studied \cite{campolongo2007effective,morris1991factorial}, and the relationship between these \emph{elementary effects} and total Sobol' indices have been investigated. For example, in \cite{kucherenko2009derivative}, a theoretical link between elementary effects and total Sobol' indices is established. More recently, Razavi et al. \cite{razavi2016new,razavi2016new-1} provide a unifying framework for the analysis of sensitivities at different scales within the input domain via a variogram-based approach. In this approach, elementary effects and Sobol' indices are recovered as special cases in the limit of small and large scales respectively. By analyzing sensitivities across the entire scale in the input space, structural features of the response surface, such as smoothness and multi-modality, are revealed.

In this paper, we are interested in the analysis of sensitivity near output extrema: which input variables---under what distributions---predominantly drive the output response near its highest/lowest levels? Output extrema are regions within the output space that are of special importance in reliability engineering; systems are often designed to operate under near-optimal conditions to maximize performance metrics. \emph{Extremum sensitivity analysis} is useful for informing the construction of tailored surrogate models that aid optimization, by isolating critical variables that most influence the output near extrema. Moreover, the quantification of uncertainties at output extrema can help engineers understand the occurrence of extreme events, and how to mitigate the associated risks. The goal of this paper is to develop indices more localized than global sensitivity metrics, automatically adapted to regions of interest within the input domain corresponding to output extrema. 

In the literature, papers directly addressing extremum sensitivity analysis are scarce. Recently, there has been work on higher-order indices utilizing skewness and kurtosis, directed at the analysis of the tails of the output distribution. In Owen et al. \cite{owen2014higher}, the authors mention that a positively skewed distribution will attain more extreme positive values than a symmetric distribution with the same variance; distributions with a large kurtosis will also attain greater extremes on both sides \cite{owen2018notitle}. Geraci et al. \cite{geraci2016high-order}, who introduce polynomial chaos-based ideas for computing skewness- and kurtosis-based indices, find applications of these indices in more accurate approximation of the probability density function (PDF) of the output response, especially near output extrema. However, neither work directly addresses the problem of sensitivity analysis near output extrema. In this work, a novel type of sensitivity indices called \emph{extremum Sobol' indices} is proposed. As the name suggests, the aim is to compute Sobol' indices of functions when the input distribution is restricted to regions in the domain leading to output extremum. This conditional distribution is modelled via a rejection sampling heuristic also known as \emph{Monte Carlo filtering} (MCF). 

Central to our investigation is the tractable computation of various sensitivity metrics. It is well-known that computing sensitivity indices---particularly in high-dimensional problems---can be prohibitive owing to the excessively large number of evaluations required. While access to a function's gradients can abate this \emph{curse of dimensionality}, our focus will be on the general case where access to a function's gradients (or adjoint) is not available. Our approach is to use orthogonal polynomial least squares approximations of the uncertain function, commonly known as generalized polynomial chaos expansions (PCE) \cite{ghanem2012stochastic,xiu2002wiener--askey,wan2006multi-element}. The application of PCE for global sensitivity analysis has been relatively well-studied and prior work can be found in \cite{sudret2008global,blatman2010efficient,sudret2015computing,seshadri2017effectively} and the references therein. The use of polynomial approximations implies the assumption that we are restricting our study to functions that can be well-approximated by globally smooth polynomials. 

The remainder of this paper is structured as follows. In \Cref{sec:coefficients} we interpret the global polynomial approximation problem as that of computing its coefficients. For polynomials of a few dimensions, this can be done without assuming \emph{anisotropy}; for higher-dimensional approximations, on the other hand, anisotropy can be used to circumvent the \emph{curse of dimensionality}. More specifically, we construct polynomial approximations over subspaces---composed by a few linear combinations of the input variables. From the coefficients in the projected space, we design a polynomial surrogate for determining the coefficients in the full space, which lead to the familiar Sobol' indices. In  \Cref{sec:extremeum}, we describe and provide algorithms for computing extremum Sobol' indices using polynomials and their ridges. The methods and utility of extremum sensitivities are demonstrated through numerical examples in \Cref{sec:num_ex}.

\section{Sensitivity analysis with polynomial ridge approximations}
\label{sec:coefficients}
Prior to detailing techniques for computing coefficients of the aforementioned polynomials, it will be worthwhile to make more precise our notion of \emph{anisotropy}. Let us consider a polynomial
\begin{equation}
p=p \left( \mM^{T} {\vx} \right) \; \; \; \text{where} \; \; \; \mM \in \mathbb{R}^{d \times n} \; \; \; \text{and} \; \; \; \vx \in \mathbb{R}^{d}.
\end{equation}
For the case where $n=d$ and where $\mM = \mI$ (the identity matrix), we have a polynomial approximation in the \emph{full} $d$-dimensional space. When $n < d$ and when the $n$ columns of $\mM$ are linearly independent, they span an $n$-dimensional subspace in $\mathbb{R}^{d}$. In this case, the function $p$ is referred to as a \emph{generalized ridge function} \cite{pinkus2015ridge}. Depending on what the \emph{weights} in $\mM$ are, the polynomial will be anisotropic with the larger weights prescribing linear combinations along which the polynomial varies the most. For notational clarity, we add a subscript below $p$ to denote the dimension over which the polynomial is constructed: $p_d$ denotes a polynomial in the full space, while $p_n$ denotes a polynomial over a subspace. From here on, we will assume that $n\ll d$. 

\subsection{Computing polynomial coefficients using ridge approximations} \label{sec:ridge_intro}
Now, let us assume the existence of a smooth function $f: \mathbb{R}^{d} \rightarrow \mathbb{R}$ with input variables $\vx$. Our goal is to compute sensitivities with respect to individual or groups of input variables. Provided one has access to $N$ input/output realizations
\begin{equation}
\left( \vx_{i}, f_i \right)_{i=1}^{N},
\end{equation}
and one can express $f$ as a polynomial (albeit approximately), we can state
\begin{equation}
\vf = 
\left[\begin{array}{c}
f_{1}\\
\vdots\\
f_{N}
\end{array}\right]\approx \left[\begin{array}{c}
p_{d} \left( \vx_{1} \right) \\
\vdots\\
p_{d} \left( \vx_{N} \right) 
\end{array}\right] = \mA_{d} \boldsymbol{\alpha},
\label{equ_1}
\end{equation}
where $\mA_{d} \in \mathbb{R}^{N \times r}$ is a \emph{polynomial design matrix} in $d$ dimensions. Writing the expansion of $p_d$ in terms of its $r$ basis terms,
\begin{equation} \label{equ:PC_expansion}
p_d(\vx) = \sum_{i=1}^r \alpha_i \Psi_i(\vx),
\end{equation}
we can identify $\mA_d$ to be
\begin{equation} \label{eqn:design_matrix}
\mA_{d} = \begin{bmatrix}
\Psi_1 (\mathbf{x}_1) & \Psi_2 (\mathbf{x}_1) & \dots & \Psi_{r} (\mathbf{x}_1) \\
\Psi_1 (\mathbf{x}_2) & \Psi_2 (\mathbf{x}_2) & \dots & \Psi_{r} (\mathbf{x}_2) \\
\vdots & \vdots & \ddots & \vdots \\
\Psi_1 (\mathbf{x}_{N}) & \Psi_2 (\mathbf{x}_{N}) & \dots & \Psi_{r} (\mathbf{x}_{N})
\end{bmatrix}.
\end{equation}
In other words, $\mA_d$ is a \emph{Vandermonde}-type matrix. The number of rows $N$ corresponds to the number of model evaluations, while the number of columns $r$ corresponds to the number of unknown polynomial coefficients. Consider two scenarios for the size of $N$ vs $r$:
\begin{itemize}
\item $N \geq r$: This leads to a polynomial least squares formulation where one tries to solve 
\begin{align}
\underset{\boldsymbol{\alpha} }{\text{minimize}} & \; \; \left\Vert \mA_{d} \boldsymbol{\alpha} - \vf \right\Vert _{2}.
\end{align}
The exponential scaling of $r$ with respect to $d$ requires $N$ to be chosen to match this scaling, implying that this approach is not computationally practical for moderate- to high-dimensional problems. In light of this, new strategies for selecting the placement of sample points have been an ongoing research topic, e.g. see \cite{seshadri2019quadrature,seshadri2017effectively}. Recently, there has also been research into leveraging multi-fidelity to reduce the number of required evaluation points \cite{palar2016multi-fidelity,jakeman2019adaptive}. 
\item $N < r$: To solve a linear system with more unknowns than equations, one has to impose some form of regularization. Assuming that the solution is sparse on a polynomial basis, Blatman and Sudret \cite{blatman2011adaptive} use least angle regression \cite{efron2004least} to identify a subset of dominant terms in the basis. The assumption of compressibility leads to formulations based on compressed sensing \cite{donoho2006compressed}, where strategies including subspace pursuit \cite{diaz2018sparse} and basis pursuit denoising \cite{tang2014subsampled,hampton2015compressive} have been proposed. The latter is an $\ell_1$-minimization strategy that uses second order cone programming to solve 
\begin{align}
\begin{split}
\underset{\boldsymbol{\alpha} }{\text{minimize}} & \; \; \left\Vert \boldsymbol{\alpha}\right\Vert _{1}\\
 \text{subject to} & \; \; \left\Vert \mA_{d} \boldsymbol{\alpha} - \vf \right\Vert_2 \leq \epsilon,
 \end{split}
\end{align}
where $\epsilon$ is a small positive constant. 
\end{itemize}
So, why the need to use a global polynomial approximation? As alluded to already, accurate polynomial approximations of $f$ permit us to leverage prior work \cite{sudret2008global} on computing global variance-based sensitivities using only the coefficient values. In this paper, we will also demonstrate how polynomial approximations are useful in computing Sobol' indices tailored to output extrema.

But what if $f$, in addition to being smooth, was also anisotropic? The idea of developing global sensitivity metrics for \emph{ridge approximations} has been studied previously by Constantine and Diaz \cite{constantine2017global}. Using ideas from \emph{active subspaces} \cite{constantine2015active}, the authors define \emph{activity scores} that are based on the eigendecomposition of the average outer product of the gradient $\nabla f$ with itself in the input space. They draw comparisons between these activity scores and total Sobol' indices (see Theorem 4.1 in \cite{constantine2017global}). While their work defines a new type of sensitivity index tailored to ridge approximations, our work builds upon the idea of ridge approximations with polynomials to propose a method for calculating moment-based sensitivity indices, such as the Sobol' indices.

Consider the \emph{polynomial ridge approximation}
\begin{equation}
\vf = 
\left[\begin{array}{c}
f_{1}\\
\vdots\\
f_{N}
\end{array}\right] \approx \left[\begin{array}{c}
p_{n} \left(\mM^{T} \vx_{1} \right) \\
\vdots\\
p_{n} \left( \mM^{T} \vx_{N} \right) 
\end{array}\right] = \mA_{n} \boldsymbol{\beta},
\label{equ_2}
\end{equation}
where $\mA_{n} \in \mathbb{R}^{N \times q}$ is another \emph{polynomial design matrix}, but now in $\mathbb{R}^n$. Since typically $q \ll r$, we will assume that $N \geq q$, implying that we have at least as many model evaluations as unknown coefficients in $\boldsymbol{\beta}$. Thus, computing $\boldsymbol{\beta}$ using regular least squares should not be a problem\footnote{We are assuming that the $N$ points are selected such that $\mA_{n}$ has a relatively low condition number.}. Note that the choice of the polynomial basis in the projected space is usually not critical. Provided the dimensionality of the projected space is small, the amount of training data required for finding the subspace is often more than sufficient for regression in the projected space. If this assumption is not true, a polynomial basis 
matching the distribution in the projected space can be constructed (see Jakeman et al. \cite[sec.~5.6.2]{jakeman2019polynomial}) to mitigate losses of accuracy at this stage. In addition, in our case, sensitivity indices and moments are not quantified in the projected space, so the basis does not need to be matched to the distribution induced in the projected space (see \Cref{sec:sobol_inds}). 

The trouble, however, is that $\boldsymbol{\beta}$ are the coefficients of $p_{n}$, and we need access to the coefficients $\boldsymbol{\alpha}$ associated with $p_{d}$ in order to estimate global sensitivities. The key is to notice that $p_n(\mM^T \vx)$ is simply a \emph{rearrangement of coefficients} from $p_d(\vx)$. In other words, we can convert the coefficients of $p_n$ to those of $p_d$ in a \emph{lossless manner}. Sampling $N \geq r$ points to evaluate the polynomial design matrices $\mA_d$ and $\mA_n$, the linear system 
\begin{equation} \label{eqn:coeff_solve}
\mA_d \bal = \mA_n \bbe,
\end{equation}
can be solved to recover the full space coefficients $\bal$. Using this method, the accuracy of the estimated polynomial coefficients is now dependent on the accuracy of estimating the dimension-reducing subspace spanned by $\mM$ and the error of a ridge approximation due to variations of $f$ outside the range of $\mM$. Assuming the presence of strong anisotropy, the subspace can be computed with significantly fewer data than a full polynomial approximation using gradient-free techniques including polynomial variable projection \cite{hokanson2018data-driven} and minimum average variance estimation (MAVE) \cite{xia2002adaptive}. In the case that the input dimension is so large that $r$ is prohibitively large for solving this linear system, subset selection schemes for truncating the full space polynomial basis---such as that proposed by Blatman and Sudret \cite{blatman2011adaptive}---can be used when determining the full space coefficients. This reduces the number of coefficients to be determined and prunes the basis of unnecessary terms, but the rearrangement step is no longer strictly lossless.

\subsection{From polynomial coefficients to Sobol' indices} \label{sec:sobol_inds}
Now, assume that there is uncertainty associated with the input variables. That is, we consider the inputs to be a set of random variables, whose distribution can be described by the PDF $\rho(\vx)$ with $\vx \in \mathcal{D} \subset \mathbb{R}^{d}$ where $\mathcal{D}$ is the support of the input PDF. In this section, it is assumed that the set of inputs is independently distributed, with marginal probability densities $\rho_j(x_j)$ for $j = 1,...,d$ such that 
\begin{equation}
\rho(\vx) = \prod_{j=1}^d \rho_j(x_j).
\end{equation}
The \emph{Sobol' indices} are then proportions of the output variance $\sigma = \text{Var}[f(\mathbf{x})]$ attributed to individual or groups of input variables \cite{sobol1993sensitivity}, defined for a set of variables $S$ as
\begin{equation} \label{eqn:sobol_defn}
\sigma_S = \frac{\text{Var}[\mathbb{E}[f(\vx)\;|\; \vx_S]]}{\sigma},
\end{equation}
where $\vx_S$ contains the variables in $S$. Directly evaluating this expression for a function via Monte Carlo requires a large number of function evaluations that scales poorly with input dimension. Below, we follow the work by Sudret \cite{sudret2008global} and show that Sobol' indices can be computed via a PCE
\begin{equation} \label{equ:random_PC_expansion}
f(\vx) \approx p_d(\vx) = \sum_{i=1}^r \alpha_i \Psi_i(\vx),
\end{equation}
given a suitable orthogonal polynomial basis $(\Psi_i)$. 

Let $(\psi^{(j)}_k)_{k\in \mathbb{N}}$ be an \emph{orthogonal basis} of univariate polynomials with respect to $\rho_j(x_j)$, where $k$ indicates the degree of the polynomial. The orthogonality of the basis implies
\begin{equation}
\int \psi^{(j)}_k(x_j) \psi^{(j)}_l(x_j)~ \rho_j(x_j)~ dx_j = \delta_{kl},
\end{equation}
where $\delta_{kl}$ is the Kronecker delta which equals to 1 if $k=l$ and 0 otherwise. Now, each multivariate basis function in \cref{equ:random_PC_expansion} $\Psi_i$ is a product of univariate orthogonal basis polynomials in each dimension, i.e.
\begin{equation}
\Psi_i(\vx) = \prod_{j=1}^d \psi^{(j)}_{i_k}(x_j),
\end{equation}
for all $i$ where $i_k$ is the degree of the polynomial in the $k$-th dimension. Thus, define a \emph{multi-index} that associates each multivariate basis polynomial with the degrees of its constituent univariate basis polynomials,
\begin{equation}
\vi = (i_i,...,i_d) \in \mathbb{N}_0^d.
\end{equation}
The polynomial series \cref{equ:random_PC_expansion} can now be rewritten using multi-indices,
\begin{equation} \label{eqn:PC_expansion_multi}
p_d(\vx) = \sum_{i=1}^r \alpha_{\vi} \Psi_{\vi}(\vx), \qquad \Psi_{\vi}(\vx) = \prod_{j=1}^d \psi^{(j)}_{i_k}(x_j).
\end{equation}
The multi-index enables us to conveniently specify the polynomials included in our basis. The set of included multi-indices is called the \emph{index set}. For instance, an isotropic index set of maximum total order $p$ is specified by
\begin{equation}
\left\lbrace \vi \in \mathbb{N}_0^d \,\, \middle| \,\, \sum_{k=1}^d i_k \leq p \right\rbrace,
\end{equation}
and a tensor grid of order $p$ is specified by
\begin{equation}
\left\lbrace \vi \in \mathbb{N}_0^d \,\, \middle| \,\, i_k \leq p  \,\,\text{for all} \, \, 1\leq k \leq d\right\rbrace.
\end{equation}
Furthermore, define a function $\text{nz}(\cdot)$ that returns the positions of the non-zero indices in a multi-index. For instance,
\begin{equation}
\text{nz}((0,1,2,0,4)) = \{2,3,5\}.
\end{equation}

Leveraging the notation introduced so far, we can express the decomposition of the (approximate) total output variance $\sigma = \text{Var}[p_d(\vx)]$ via the following equation 
\begin{equation} \label{equ:anova}
\sigma = \myEx\left[(p_d(\vx) - \myEx \left[p_d(\vx)\right])^2\right] = \int_{\mathcal{D}} \left(\sum_{i=2}^r \alpha_{\vi} \Psi_{\vi}(\vx)\right)^2 \rho(\vx)~ d\vx.
\end{equation}
Owing to orthonormality, the polynomial series expression \cref{eqn:PC_expansion_multi} allows us to compute the Sobol' indices straightforwardly given the polynomial coefficients $\bal$ by comparison with the Sobol' decomposition
\begin{equation} \label{eqn:sobol_decomp}
f(\vx) = \sum_{S \subseteq [d]} f_S(\vx_S)
\end{equation}
where $[d] = \{1,2,...,d\}$, $\vx_S$ selects all coordinates of $\vx$ belonging to the subset of indices $S$, and $f_S(\vx_S)$ is recursively defined as
\begin{equation}
f_S(\vx_S) = \mathbb{E} [f(\vx) \; | \; \vx_S] - \sum_{T \subset S} f_T(\vx_T).
\end{equation}
Using this decomposition, Sudret shows that the Sobol' index for a set of variables indexed by $S$ is given by
\begin{equation} \label{equ:sobol}
\sigma_S = \frac{\sum_{\text{nz}(\vi) = S} \alpha_{\vi}^2}{\sum_{\vi} \alpha_{\vi}^2}.
\end{equation}
A measure of the importance of an input variable is the total Sobol' index \cite{sobol2001global}, which is the sum of all Sobol' indices that contain the variable concerned. This index factors in all orders of interaction within the function that involves the concerned variable. For instance, for variable $x_i$, the corresponding total Sobol' index is
\begin{equation} \label{eqn:total_sobol}
\hat{\sigma_i} = \sum_{i\in S} \sigma_S.
\end{equation}

\section{Extremum sensitivity analysis}
\label{sec:extremeum}
In the previous section, we described a method to use coefficients obtained from a polynomial ridge approximation to calculate Sobol' indices. In the following section, the second key contribution of this paper---the proposal of polynomial-based sensitivity indices focusing on inference near output extrema---will be addressed.

\subsection{Skewness-based indices} \label{sec:skew_indices}
First, we review two works in literature related to global sensitivity analysis based on the decomposition of skewness.  Assuming that the input domain is a unit hypercube $\mathcal{D} = [0,1]^d$ endowed with uniform marginal distributions, Owen et al. \cite{owen2014higher} define \emph{higher-order Sobol' indices} by generalizing the Sobol' identity, which is based on the functional ANOVA (Sobol') decomposition \eqref{eqn:sobol_decomp}. This identity states that \cite[Thm.~2]{sobol2001global}
\begin{equation}
\sigma\sigma_S = \left(\int_{\mathcal{D}} \int_{\mathcal{D}} f(\vx) f(\vx_S;\vx^{(1)}_{-S}) \; d\vx \;d\vx^{(1)} - \mu^2 \right),
\end{equation}
where the notation $\vx_S;\vx^{(1)}_{-S}$ denotes a point $\vy$ in $\mathcal{D}$ where $y_j = x_j$ for $j \in S$ and $y_j = x^{(1)}_j$ for $j \notin S$, and $\mu$ denotes the mean of the output. The generalization for calculating the (unnormalized) skewness index for a subset of variables $S$ reads
\begin{equation}
t_S' = \int_{\mathcal{D}}\int_{\mathcal{D}}\int_{\mathcal{D}} f(\vx) f(\vx_S;\vx^{(1)}_{-S}) f(\vx_S;\vx^{(2)}_{-S}) \; d\vx \;d\vx^{(1)} \;d\vx^{(2)} - \mu^3.
\end{equation}
Also in \cite{owen2014higher}, it is shown that 
\begin{equation}
t_S' + \mu^3 = \mathbb{E} \left[\mathbb{E} [f(\vx)~|~ \vx_S]^3\right],
\end{equation}
which parallels the definition of Sobol' indices where 
\begin{equation}
\sigma \sigma_S + \mu^2 = \mathbb{E} \left[\mathbb{E} [f(\vx) ~|~ \vx_S]^2\right].
\end{equation}
The authors in \cite{owen2014higher} remark that skewness indices may provide indication that certain variables are important for attaining extreme values of the output, with the \emph{sign} of the skewness indices indicating whether the variables contribute significantly to output maximum (positive skewness index) or minimum (negative skewness index). However, they do not provide further theoretical justification for this claim. 

The formulation described by Owen et al. lends straightforwardly to computation by Monte Carlo samples, but the high dimensionality of the integrals implies that many samples are required. In Geraci et al. \cite{geraci2016high-order}, skewness indices are defined by evaluating the third central moment of the Sobol' decomposition \eqref{eqn:sobol_decomp} and expanding the sum, resulting in indices that are different from those defined in \cite{owen2014higher}. For brevity, we omit the full exposition and refer the reader to their work for further details. The authors provide an approach for evaluating these indices via polynomial chaos, by expressing the indices in terms of polynomial coefficients. However, unlike the case for variance-based indices, the resultant indices cannot be expressed in a form as simple as \cref{equ:sobol}. Here, we cite the full result from \cite{geraci2016high-order} for the skewness index of a set of variables indexed by $S$. Letting $\vp, \vq, \vr$ be the corresponding multi-indices to the scalar indices $p, q, r$ respectively, we define the following sets
\begin{equation*}
S^1 = \{p \, : \, \text{nz}(\vp) = S \},
\end{equation*}
\begin{equation*}
S^2 = \{(p,q) \, : \, \text{nz}(\vp) \cup \text{nz}(\vq) = S\},
\end{equation*}
\begin{equation*}
S^3 = \{(p,q,r) \, : \, \text{nz}(\vp) \cup \text{nz}(\vq) \cup \text{nz}(\vr) = S \}.
\end{equation*}
It can be shown that the skewness index for $S$ is
\begin{align} \label{eqn:skewness_indices}
\begin{split}
t_S = \frac{1}{\sigma^{3/2} \gamma} &\left( \sum_{p \in S^1} \alpha_\mathbf{p}^3 \myEx\left[ \Psi_\mathbf{p}^3 \right] +3 \sum_{\substack{(p,q) \in S^2 \\ p\neq q}} \alpha_\mathbf{p}^2 \alpha_\mathbf{q} \myEx\left[ \Psi_\mathbf{p}^2 \Psi_\mathbf{q}\right] + 6 \sum_{\substack{(p,q,r) \in S^3 \\ q>p}} \alpha_\mathbf{p} \alpha_\mathbf{q} \alpha_\mathbf{r} \myEx\left[ \Psi_\mathbf{p} \Psi_\mathbf{q} \Psi_\mathbf{r}\right]\right).
\end{split}
\end{align}
The total skewness index for $x_i$ is defined analogously to \cref{eqn:total_sobol}
\begin{equation}
\hat{t}_i = \sum_{i\in S} t_S.
\end{equation}
In the same work, the authors also establish the importance of these indices in modelling the tails of the output PDF, but make no explicit reference to sensitivity analysis near output extrema.

\subsection{Extremum Sobol' indices} \label{sec:mcf_indices}
Although skewness indices pertain to the analysis of the tails of the output distribution, limitations exist with respect to the use of these indices for extremum sensitivity analysis. For a simple example, consider the following function
\begin{equation}
f(x_1, x_2) = x_1^2 + 100x_2, \qquad (x_1, x_2) \sim \mathcal{U}\left[-1, 1\right]^2.
\end{equation}
By inspection, it is clear that $x_2$ is responsible for none of the (small but non-zero) output skewness, but it is more important than $x_1$ throughout the entire input domain, including regions of output extrema. 

In this paper, we propose a novel set of variance-based metrics that quantify the relative input sensitivities near regions in the domain yielding extreme output values, called \emph{extremum Sobol' indices}. The aim is to analyze the variance decomposition of the function conditioned near its extrema, producing indices that have a clear interpretation for extremum sensitivity analysis. To characterize the input distribution corresponding to output conditioning, a heuristic related to Monte Carlo filtering (MCF) \cite[p.~39]{saltelli2008global} is used. In MCF, the input space is partitioned into two regions---one where the corresponding outputs fall in a range of interest $B$, and one where this does not happen, $\bar{B}$. The importance of a variable $x_i$ is determined by the difference between the distributions of $x_i$ conditioned within $B$ and its complement $\bar{B}$ respectively, where the distributions are inferred using Monte Carlo. In our method, MCF is used to isolate regions containing points leading to output maximum and minimum, and the Sobol' indices of the quantity of interest are evaluated under the conditional distributions restricted to these regions. Sobol' indices resulting from points leading to output maximum are called \emph{top Sobol' indices}, and those from points leading to output minimum are called \emph{bottom Sobol' indices}.

Computation of extremum Sobol' indices involves the following steps. 
\begin{enumerate}
\item \textbf{Global polynomial approximation}: Fit a global polynomial approximation to the model by input/output pairs $\{\vx, f(\vx)\}$ via sampling from $\mathcal{D}$ according to $\rho(\vx)$. Note that this global approximation can take the form of a ridge approximation.
\item  \textbf{Rejection sampling}: Generate another $N$ random input points in $\mathcal{D}$ according to $\rho(\vx)$. Estimate the maxima and minima by pushing these samples through the polynomial approximation. Isolate input-polynomial output $\left\{\vx, p(\vx) \right\}$ pairs that fall within the top 5\% and bottom 5\% of the maxima and minima respectively.
\item  \textbf{Marginals and correlation determination}: Approximate the marginal distributions for each parameter for the two groups using input/output pairs. Fit a tailored bandwidth kernel density estimate to these marginals. Additionally, evaluate the empirical correlation matrix for these two groups.
\item  \textbf{Extremum Sobol’ indices}: Construct polynomial chaos approximations over the groups, using additional input/output pairs $\{\vx, f(\vx)\}$ sampled from the correlated distribution over the respective groups. Then, use \Cref{alg:extr_sobols} to compute the extremum Sobol’ indices for each group.
\end{enumerate}

For the final step, it is important to note that the extremum points resulting from MCF are in general distributed with a correlated measure, even if the original input measure is independent. In this case, the output variance cannot be decomposed in a similar manner to the Sobol' decomposition, and the indices computed from the definition \eqref{eqn:sobol_defn} mixes the contribution of variables outside of the variables in $S$; the sum of all Sobol' indices thus does not result in unity. Therefore, we compute sensitivities using the formulation proposed by Li et al. \cite{li2010global} which involves a covariance decomposition, for which Navarro et al. \cite{navarro2014polynomial} expanded upon with an algorithm based on polynomial chaos expansions. \Cref{alg:extr_sobols} in \ref{sec:extr_sobol} summarizes the steps of this computation given extreme points isolated by MCF.

To calculate the extremum Sobol' indices using a polynomial chaos approximation, a polynomial basis orthogonal to the correlated measure determined in step 4 above is required. There are different methods to construct orthogonal polynomial bases in correlated spaces. For instance, one can map the dependent variables to an independent space via the Nataf transform \cite{liu1986multivariate}. In \Cref{alg:extr_sobols}, we use the Gram-Schmidt method \cite{witteveen2006modeling,jakeman2019polynomial} to generate evaluations on a polynomial basis orthogonal to the correlated measure. We start with the polynomial design matrix $\mA_u$ as in equation \eqref{eqn:design_matrix}, with the $N$ points of evaluation $\vx_i$ sampled from the correlated measure. The basis polynomials here can be taken as ones orthogonal to the full input space measure $\rho(\vx)$.  The sampling step can be performed via the Nataf transform. Then, the QR decomposition
\begin{equation}
\frac{1}{\sqrt{N}} \mA_u = \mQ_u \mR
\end{equation}
is calculated, where $\mQ_u$ has orthogonal columns, and $\mR$ is upper triangular. Because of the distribution of the sample points, this step is a Monte Carlo approximation to the orthogonalization of the polynomial basis with respect to the inner product
\begin{equation}
\langle \Psi_1, \Psi_2 \rangle = \int \Psi_1(\vx) \Psi_2(\vx) ~\omega(\vx) ~ d\vx
\end{equation}
where $\omega$ is the correlated (extremum) measure and $\Psi_1$ and $\Psi_2$ are two polynomial basis functions. Polynomial coefficients can be solved on this basis using least squares and Sobol' indices quantified via the coefficients. However, note that the Gram-Schmidt process produces basis polynomials that mix contributions from different variables. Therefore, the Sobol' indices need to be calculated via Monte Carlo on the basis polynomials instead of being simply read out from \eqref{equ:sobol}. For this algorithm, we remark that additional function evaluations are required only for solving the coefficients in the new polynomial basis; all Monte Carlo steps are performed on evaluations of the polynomial basis terms.

\section{Numerical examples}
\label{sec:num_ex}

In the following section, polynomial ridge approximations and extremum sensitivity analysis are applied on several numerical case studies. 
We aim to highlight the following points:
\begin{enumerate}
\item Demonstrate the estimation of Sobol' indices using polynomial ridge functions via the least squares formulation described in \Cref{sec:coefficients}.
\item Discuss empirical parallels between skewness indices (\Cref{sec:skew_indices}) and extremum Sobol' indices (\Cref{sec:mcf_indices}) at capturing important features at top and bottom output values.
\item Compare the computational cost of our framework with Monte Carlo and other polynomial chaos methods such as adaptive LARS where applicable.
\end{enumerate}

\subsection{Ridge approximation of an analytical function}
In this example, we study the estimation of Sobol' indices of the function $f: [-1,1]^6 \rightarrow \myRe$, where
\begin{equation} \label{eqn:f_def}
f(\vx) = f_1(\mU^T \vx) + s f_2(\mW^T \vx),
\end{equation}
with
\begin{equation}
f_1(\vu) = u_1^2 - 0.1u_2^2 - 2u_1^2u_2^2 + 2u_1u_2, \qquad \vu \in \myRe^2,
\end{equation}
\begin{equation}
f_2(\vw) = w_1^2 - 0.1w_2^2 - 2w_1^2w_2^2 + 2w_1w_2, \qquad \vw \in \myRe^4,
\end{equation}
and the orthogonal columns of $\mU, \mW$ span, respectively, the column and left null space of
\begin{equation}
\begin{bmatrix}
2 & 3 & 1 & -4 & 1 & 0.1 \\
-3 & 1 & -2 & 1 & -1 & -0.3 
\end{bmatrix}^T.
\end{equation}
Here, $s$ is a parameter that modulates the deviation of $f(\vx)$ from an exact ridge function within the column space of $\mU$---i.e. it allows the function to vary slightly outside of the plane spanned by the columns of $\mU$. We evaluate the sensitivity indices from noisy samples $\{\vx_i, f(\vx_i) + \epsilon_i\}_{i=1}^N$, where the sample inputs $\vx_i$ are drawn from the uniform distribution over the input domain, $\rho(\vx)= 1/2^6$, and $\epsilon_i$ represents additive Gaussian random noise with a standard deviation of $10^{-4}$. Five methods are compared for this task:
\begin{itemize}
\item Polynomial ridge approximation: Using MAVE \cite{xia2002adaptive}---a gradient-free dimension reduction method---we estimate a dimension-reducing subspace for the function of interest. Then, we fit a low-dimensional quadratic polynomial ridge over the reduced space $p_n$ with $n=2$, whose coefficients we use to calculate the coefficients in the full space via \cref{eqn:coeff_solve}. Since $n$ is small, we can afford to use a tensor grid. A Legendre polynomial basis is used in the projected space; since we do not need to quantify moments, the choice of basis is not critical here. Note that in this example, the choice of $n=2$ is prescribed according to the known function, but in practice for an unknown model, methods such as cross validation (or visualization with projection plots) may be used for determining the suitable dimensionality.
\item Quasi-Monte Carlo: We compute the integrals defining the Sobol' indices via a Sobol' sequence, as specified in \cite{sobol2001global}. We adapt the open-source Python library \texttt{SALib} \cite{herman2017salib} for this method.
\item Adaptive least angle regression (LARS): Based on an approach by Blatman and Sudret \cite{blatman2011adaptive} (see \ref{sec:appenB}), we evaluate the Sobol' indices via orthogonal polynomial regression over a pruned basis, which is a subset of an isotropic total order basis of maximum degree 4 over the full space. The open-source Python library \texttt{scikit-learn} \cite{pedregosa2011scikit-learn} is used for LARS.
\item Fitting a polynomial on the full input space via least squares with a total order basis of maximum order 3.
\item Fitting a polynomial on the full input space via least squares with a total order basis of maximum order 4.  
\end{itemize}

\begin{figure} 
\begin{minipage}{0.5\textwidth}
\centering
\includegraphics[width=.95\textwidth]{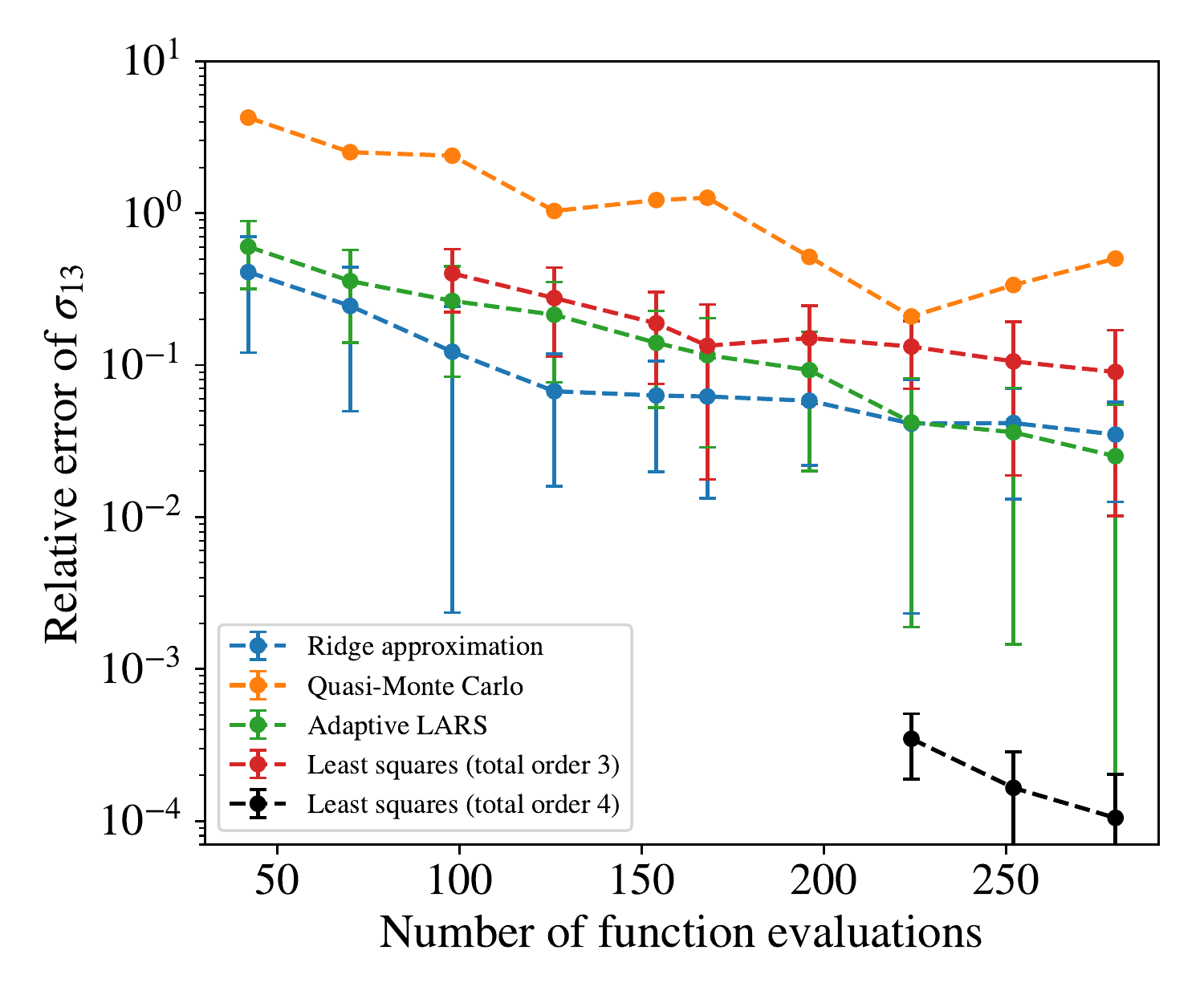}
\end{minipage}
\hfill
\begin{minipage}{0.5\textwidth}
\centering
\includegraphics[width=.95\linewidth]{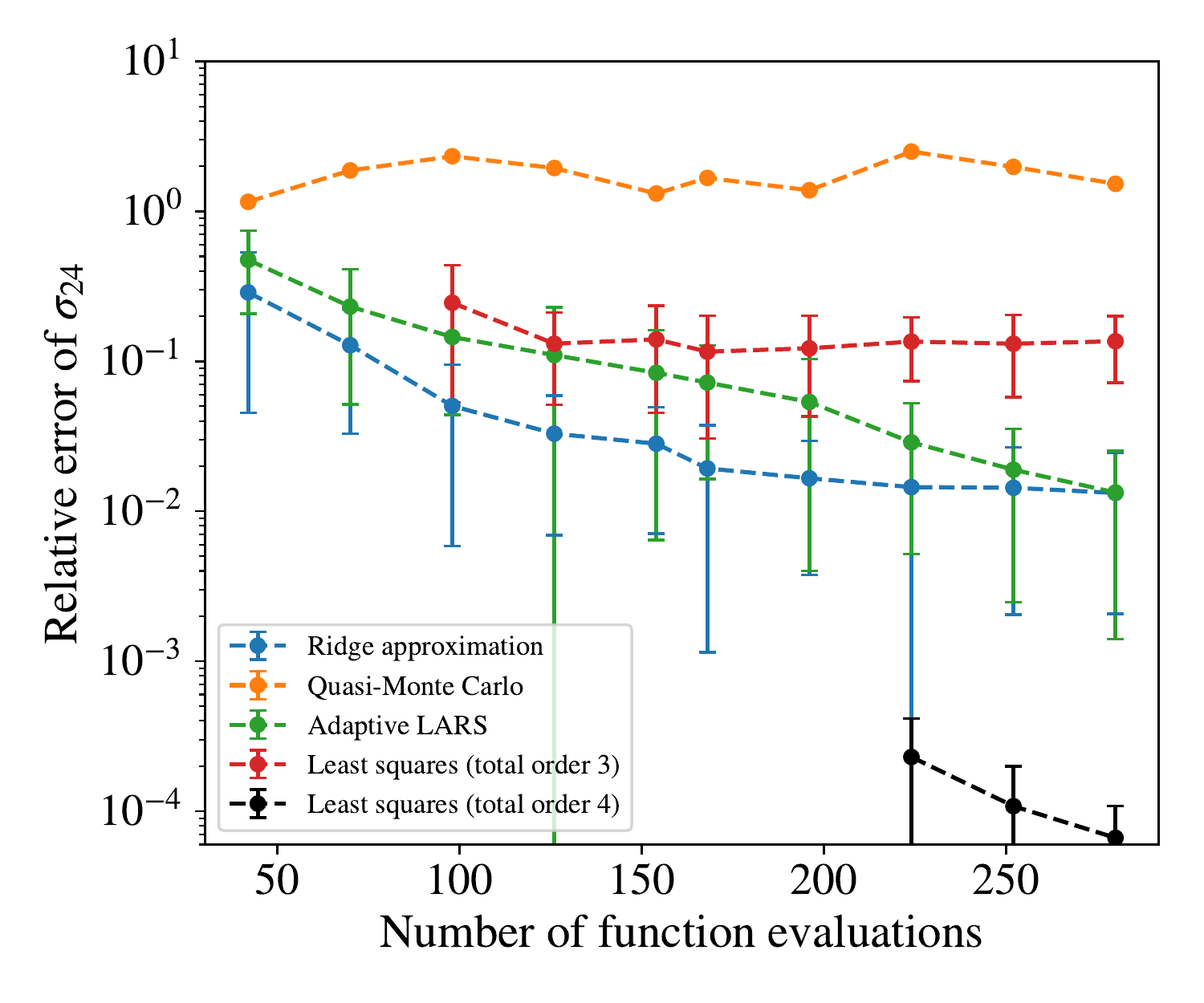}
\end{minipage}
\begin{minipage}{0.5\textwidth}
\centering
\includegraphics[width=.95\textwidth]{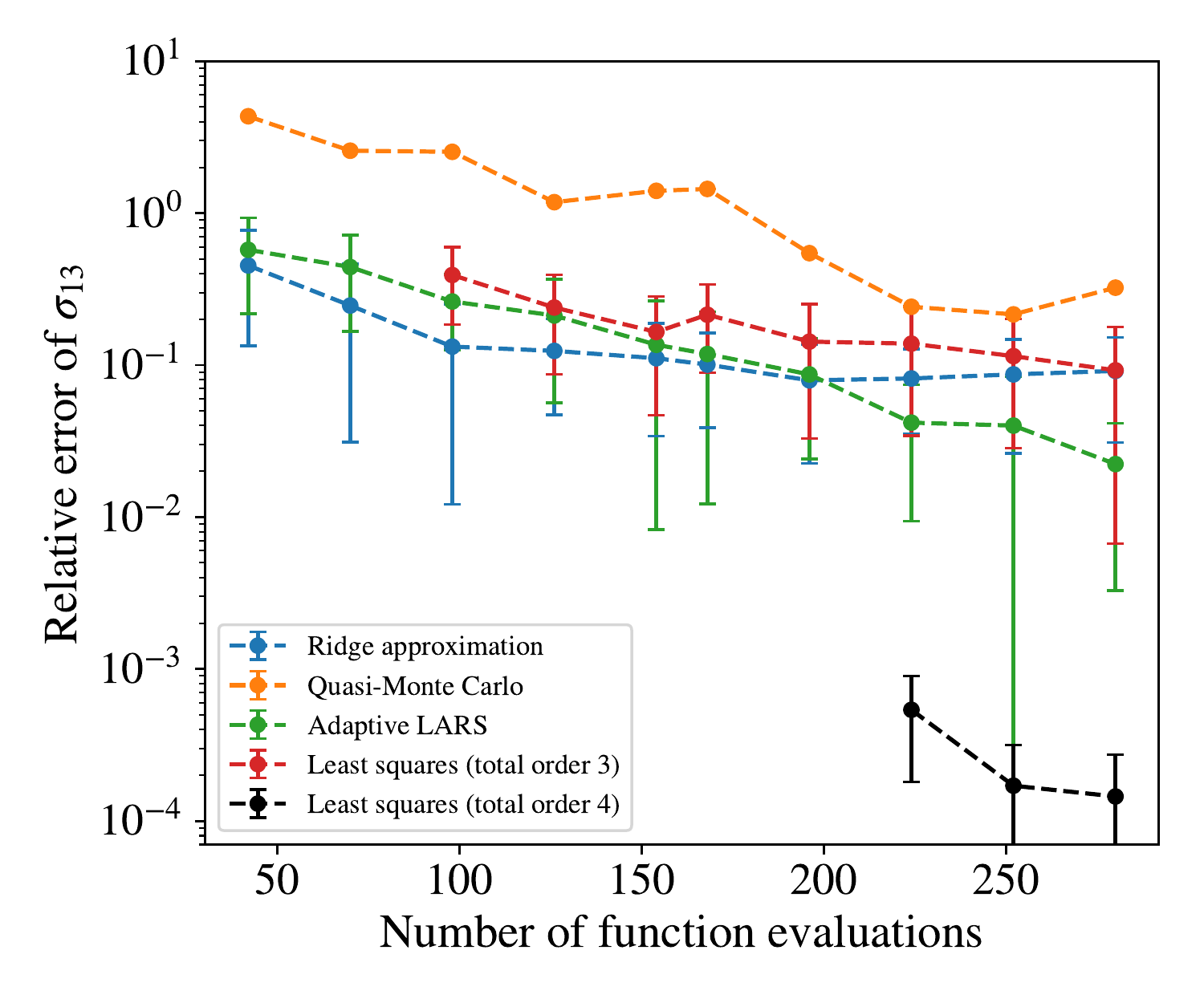}
\end{minipage}
\hfill
\begin{minipage}{0.5\textwidth}
\centering
\includegraphics[width=.95\linewidth]{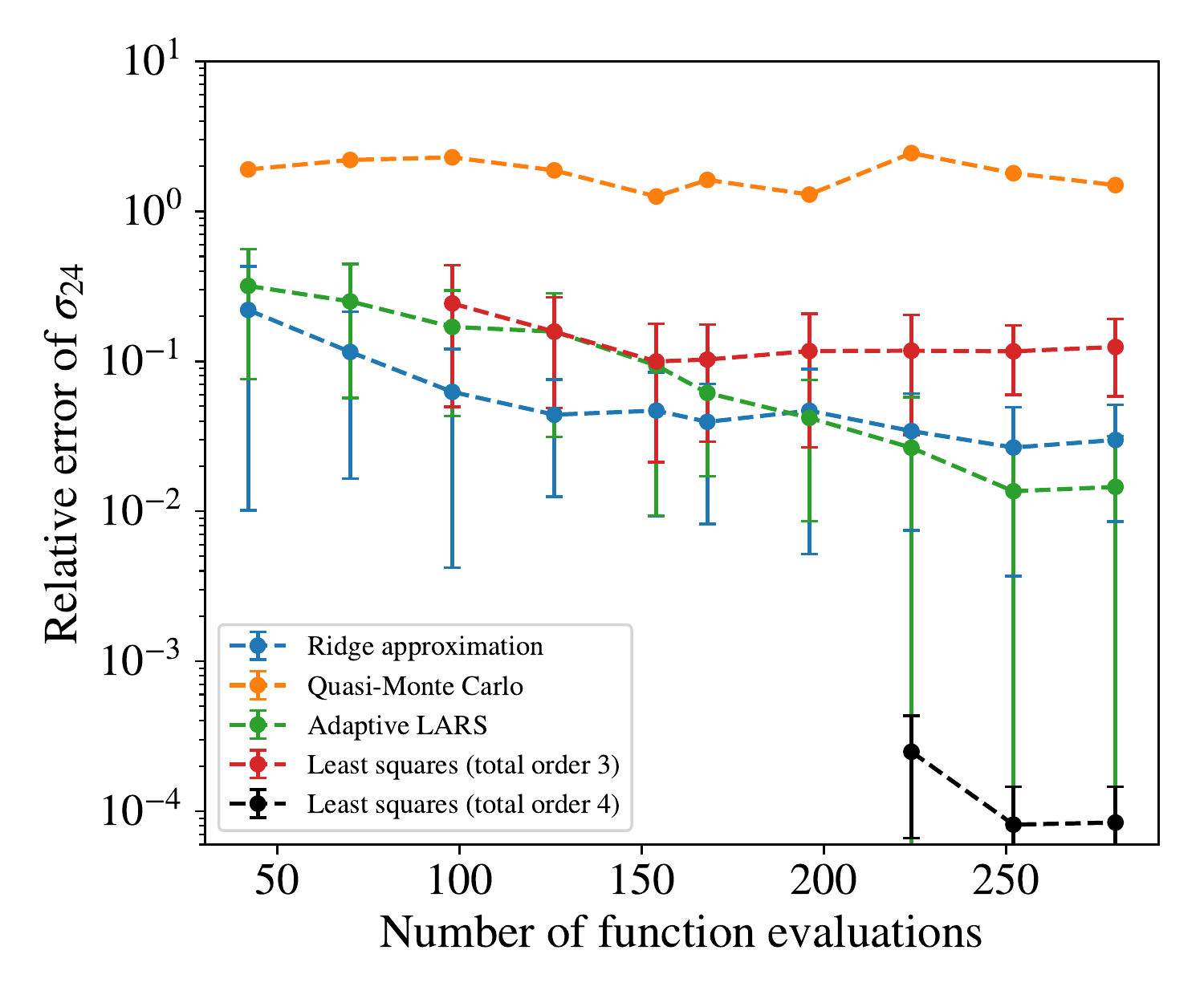}
\end{minipage}

\caption{Accuracy of the second-order Sobol' indices $\sigma_{13}$ (left) and $\sigma_{24}$ (right) for the function $f(\vx)$ defined in \cref{eqn:f_def}, averaged over 30 trials. The parameter $s$ is set to be 0 (top) and 0.1 (bottom) respectively.}
\label{fig:ridge_results}
\end{figure}

Note that the final two approaches place a lower limit on the number of samples required (84 and 210 respectively) as explained in \cref{sec:ridge_intro}. The two largest Sobol' indices $\sigma_{13}$ and $\sigma_{24}$ are computed using the methods above. To benchmark these results, we also compute the quantities using integration via Gauss quadrature on a degree 4 tensor grid of Legendre polynomials with $(4+1)^6 = 15625$ points, the results of which we treat as the truth. The fact that $f(\vx)$ is a polynomial allows this integration to be exact. \Cref{fig:ridge_results} shows the results of this comparison for $s=0$ and $s=0.1$, respectively. 

From the results, it is clear that polynomial-based techniques offer a significant advantage over quasi-Monte Carlo integration. This shows that sensitivity analysis of a function that is smooth and well approximated by a polynomial can be greatly facilitated by our approach. Additionally, the fact that this function is (approximately) a ridge function gives the ridge approximation approach an advantage over other sparse methods, including adaptive least angle regression. In general, ridge functions need not have sparse coefficients unless the ridge structure aligns well with the canonical directions in the input space. The case of $s=0.1$ adds variation outside of the two-dimensional subspace spanned by the columns of $\mU$, which reduces the accuracy of the ridge approximation. Nonetheless, the ridge approximation approach fares well against other sparse approaches. However, if we are allowed to evaluate the function more than 210 times, regression on the full space using a total order grid is superior. 

\subsection{Extremum sensitivity analysis of an analytical function} \label{sec:ex_analytical}
\begin{figure}
\centering
\includegraphics[width=1\textwidth]{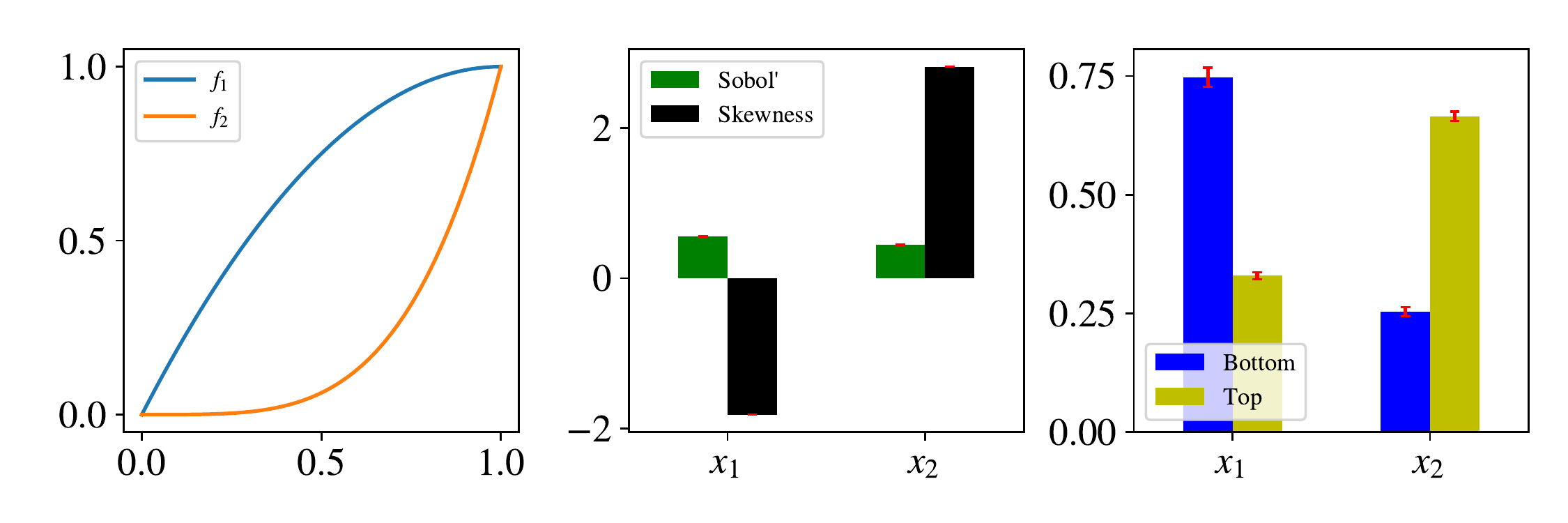}
\caption{The plots for the two subfunctions $f_1$ and $f_2$ in \eqref{eqn:extr_fn} (left), and the associated total Sobol' and skewness (middle) and extremum Sobol' indices (right). }
\label{fig:analytical_extr}
\end{figure}

Consider the additive function $f: [0, 1]^2 \rightarrow \mathbb{R}$ where
\begin{equation} \label{eqn:extr_fn}
f(\vx) = f_1(x_1) + f_2(x_2) = -x_1(x_1 - 2) + x_2^4, \qquad \vx \sim \mathcal{U}[0,1]^2.
\end{equation}
We calculate total sensitivity indices for this function via polynomial methods, which are exact in this case because the function is a polynomial\footnote{However, Monte Carlo is still used to calculate the extremum Sobol' indices after polynomials are fitted near output extrema. Sufficient samples are used to ensure that the uncertainty is small.}. Sobol' indices, skewness indices and extremum Sobol' indices are plotted for this function in \Cref{fig:analytical_extr}. In all bar chart figures, the bar heights show the mean and red error bars show $\pm1$ standard deviation over 30 trials. It can be seen from the middle plot that both variables have similar Sobol' indices, but there is a clear difference for the bottom and top Sobol' indices. For low output values, $f_2$ is nearly flat, implying that it is not influential to the function value; at high outputs, $f_1$ is relatively flat, implying that $f_2$ should be the more important variable. These pieces of information are revealed from the extremum Sobol' indices. 

From the skewness indices, it can be seen that $x_1$ is responsible for skewing the output distribution to the right, and $x_2$ to the left. In this example, it can be observed that $x_1$, with a negative skewness index, is important near the output minimum; for $x_2$, it is important near the output maximum, and has a positive skewness index. This agrees with the interpretation of skewness indices by Owen et al. \cite{owen2014higher} mentioned in \Cref{sec:skew_indices}, although it should be noted that their definition of skewness indices is different from the ones computed here. However, whether the sign of skewness indices are indicators of extremum sensitivities in general, and its precise connection with extremum Sobol' indices remain open topics.

\subsection{Borehole function}

In this case study, we apply methods of extremum sensitivity analysis on the borehole function, which models the steady-state water flow rate for a hypothetical borehole geometry \cite{morris1993bayesian}. Our objective here is to identify which variable is responsible for the greatest variation over extreme values of the steady-state water flow. The model is specified by eight input parameters assumed to be independent. Their distributions are given in \Cref{tab:borehole_inputs}. The water flow rate depends on the parameters as
\begin{equation}
f(\vx) = \frac{2\pi T_u (H_u - H_l)}{\ln\left(\frac{r}{r_w}\right) \left(1+ \frac{2LT_u}{\ln(r/r_w) r_w^2 K_w} + \frac{T_u}{T_l}\right)}.
\end{equation}
\begin{table} []
\centering
\begin{tabular}{@{}cccc@{}}
\hline
\textbf{Input} & \textbf{Distribution}&\textbf{Distribution parameters} & \textbf{Description} \\
\midrule
$r_w$ & Gaussian & 0.10, 0.0161812 & Radius of borehole (m) \\

$r$ & Uniform& 100, 50000 & Radius of influence (m) \\

$T_u$ & Uniform& 63070, 115600 & Transmissivity of upper aquifer  (m$^2$/yr) \\

$H_u$& Uniform & 990, 1110 & Potentiometric head of upper aquifer (m) \\

$T_l$& Uniform & 63.1, 116 & Transmissivity of lower aquifer (m$^2$/yr) \\

$H_l$ & Uniform& 700, 820 & Potentiometric head of lower aquifer (m)\\

$L$& Uniform & 1120, 1680 & Length of borehole (m)\\

$K_w$ & Uniform& 1500, 15000 & Hydraulic conductivity of borehole (m/yr)\\
\bottomrule
\end{tabular}
\caption{Input variables to the borehole function. Distribution parameters to a Gaussian distributed variable refer to the mean and standard deviation respectively; for a uniformly distributed variable, they refer to the lower and upper limits of the support.}
\label{tab:borehole_inputs}
\end{table}

\begin{figure}
\centering
\includegraphics[width=1\textwidth]{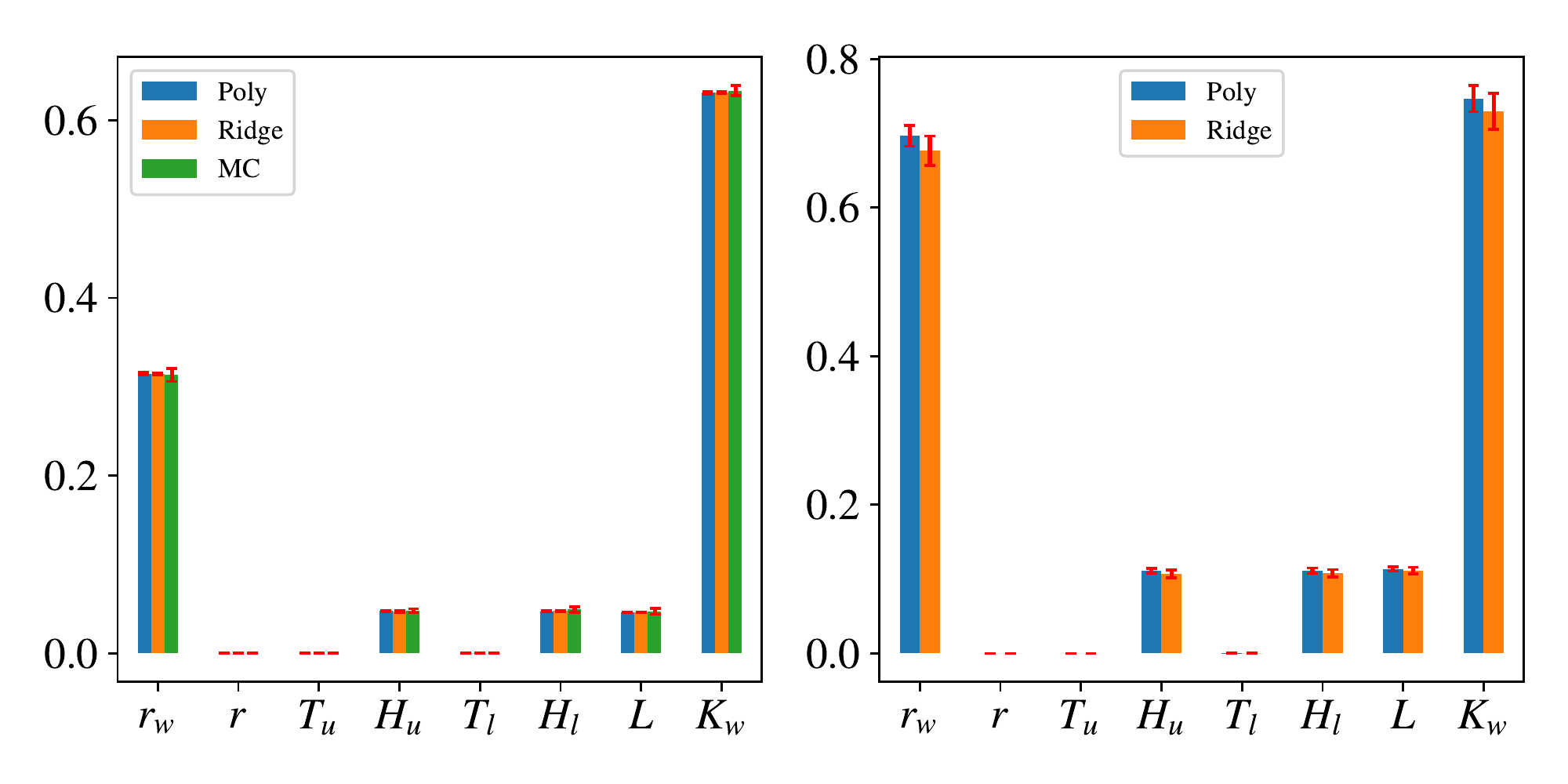}
\caption{Total Sobol' (left) and skewness (right) indices for the borehole function. }
\label{fig:borehole_inds}
\end{figure}

Treating the water flow rate as the quantity of interest, we compute the total Sobol' and skewness indices with respect to the input variables. Three methods of computation are compared:
\begin{enumerate}
\item (Poly) Using an orthogonal polynomial least squares approximation on an isotropic basis with maximum total order of 3 (with a cardinality of 165), trained from 300 input/output pairs.
\item (Ridge) Using a polynomial ridge approximation of degree 3 in a 4-dimensional subspace\footnote{In this example and the next, the subspace dimension is found by trial-and-error, and the lowest value that results in acceptable performance is selected. In practice, methods such as cross validation (or visualization with projection plots) can need to be used for determining the suitable dimensionality, but we omit a detailed study here for brevity.}, whose coefficients are subsequently converted for a polynomial on an isotropic basis with maximum total order of 3. The subspace is found via polynomial variable projection \cite{hokanson2018data-driven} using 150 input/output pairs. In this and the next numerical examples, Legendre bases on tensor grids are used for fitting polynomials in projected spaces. 
\item (MC) Using $10^5$ model evaluations to calculate the indices via the Monte Carlo method by Sobol' \cite[Thm. 2]{sobol2001global}.
\end{enumerate}
The training set $\{\vx_i, f(\vx_i)\}_{i=1}^N$ is formed by sampling the inputs independently according to the corresponding marginal input distributions, and evaluating the corresponding output values. In \cref{fig:borehole_inds}, results for the computations are shown. The third method is omitted for the skewness indices because no equivalent Monte Carlo formulation is proposed for skewness indices. The algorithm for computing skewness indices via polynomial chaos is explained in \ref{sec:appen}. From this figure, it can be seen that the Sobol' indices computed using a polynomial approximation matches well with the Monte Carlo estimates, and using a ridge approximation in lieu of a full space polynomial reduces the amount of training data required while having minor effects on the accuracy of the computed indices. The choice of polynomial degrees for this and the next example is confirmed by verification results found in \ref{sec:ver}.

Next, we compute the extremum Sobol' indices using the procedure described in \Cref{sec:mcf_indices}. The total top and bottom Sobol' indices are calculated using MCF to select 5\% of the input points yielding maximum and minimum outputs respectively out of $10^5$ points. Three variants of the computation algorithm are tested differing in the form of the model used to generate the pool of samples for MCF at step 2---namely 
\begin{enumerate}
\item Using $10^5$ function evaluations directly (instead of a polynomial approximation as specified in step 1);
\item Using a global polynomial model at total degree 4 (cardinality 495) fitted with 700 observations; and
\item Using a ridge polynomial at total degree 4 in a 4-dimensional subspace fitted with 400 observations.
\end{enumerate}
Note that the polynomial degree is increased by one compared to the approximation models used to calculate global sensitivity indices previously in this section. This is to ensure better accuracy at the tails of the output. For all methods, an additional 1400 function evaluations are used to fit a degree 5 polynomial at the extrema.

\begin{figure}
\centering
\includegraphics[width=1\textwidth]{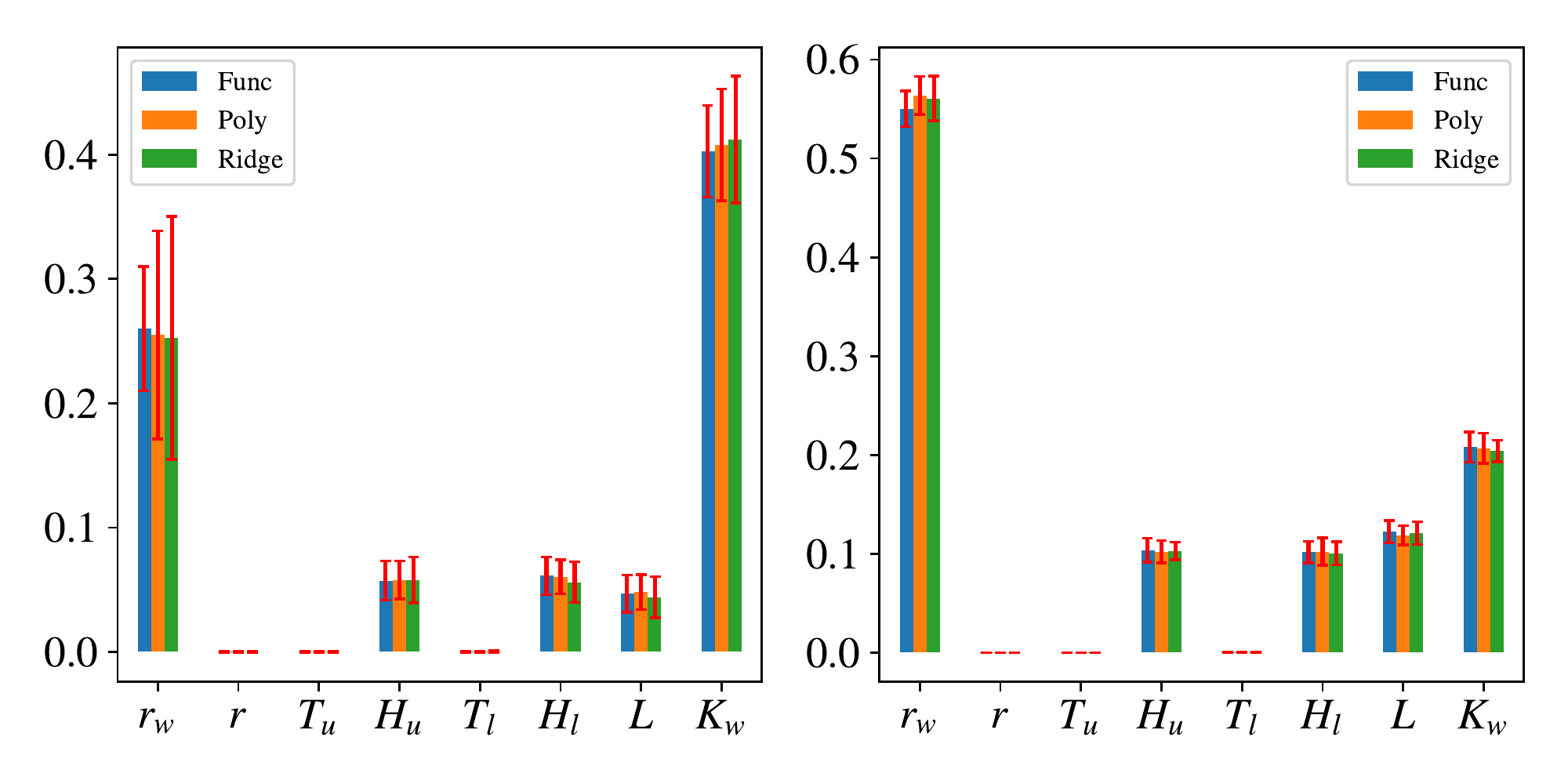}
\caption{Bottom (left) and top (right) total Sobol' indices for the borehole function using function evaluations (Func), a polynomial approximation (Poly) and a ridge approximation (Ridge).}
\label{fig:borehole_extr}
\end{figure}

\begin{figure}
\centering
\includegraphics[width=1\textwidth]{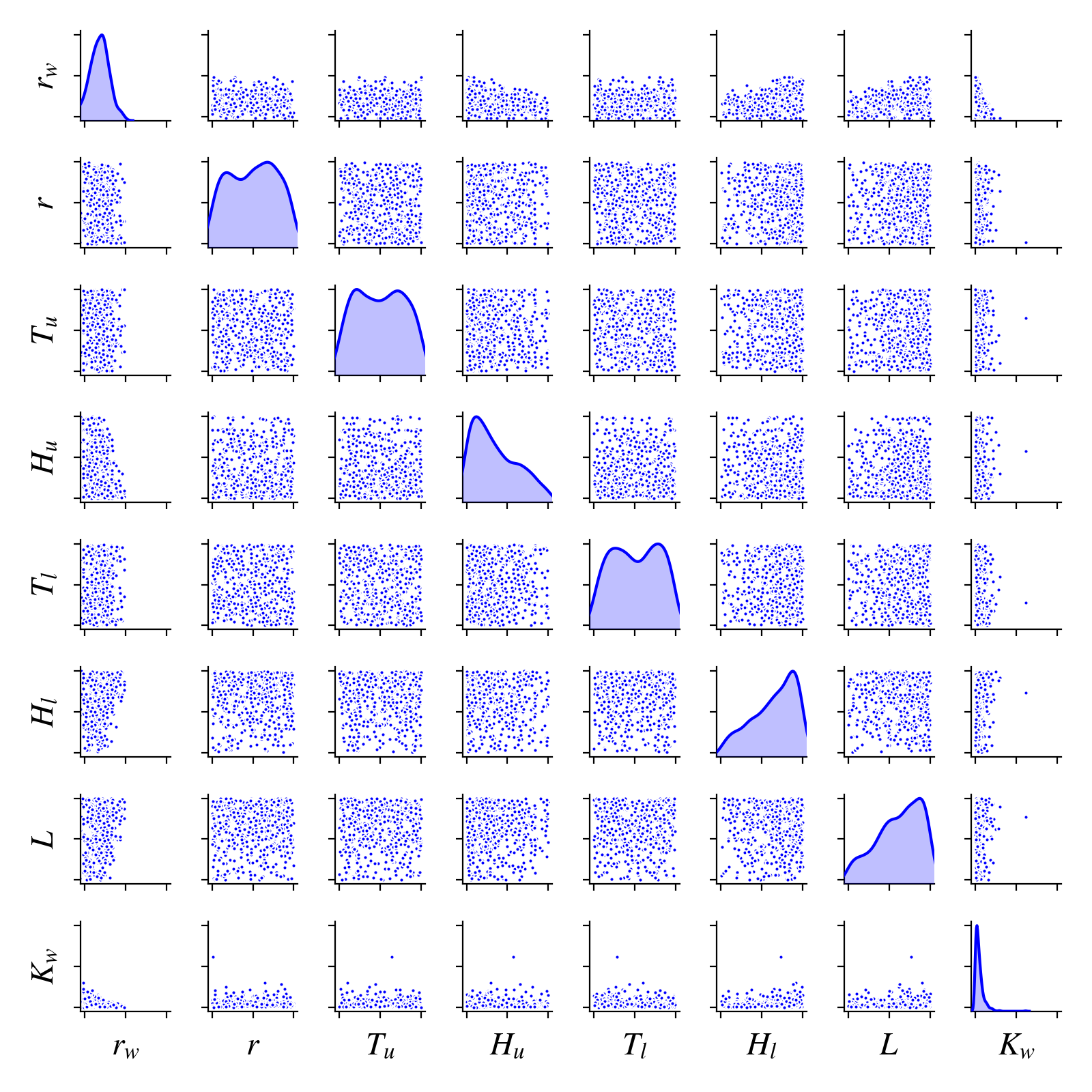}
\caption{Pairwise scatter plot of isolated samples leading to low outputs for the borehole model.}
\label{fig:borehole_bot_pdf}
\end{figure}

\begin{figure}
\includegraphics[width=1\textwidth]{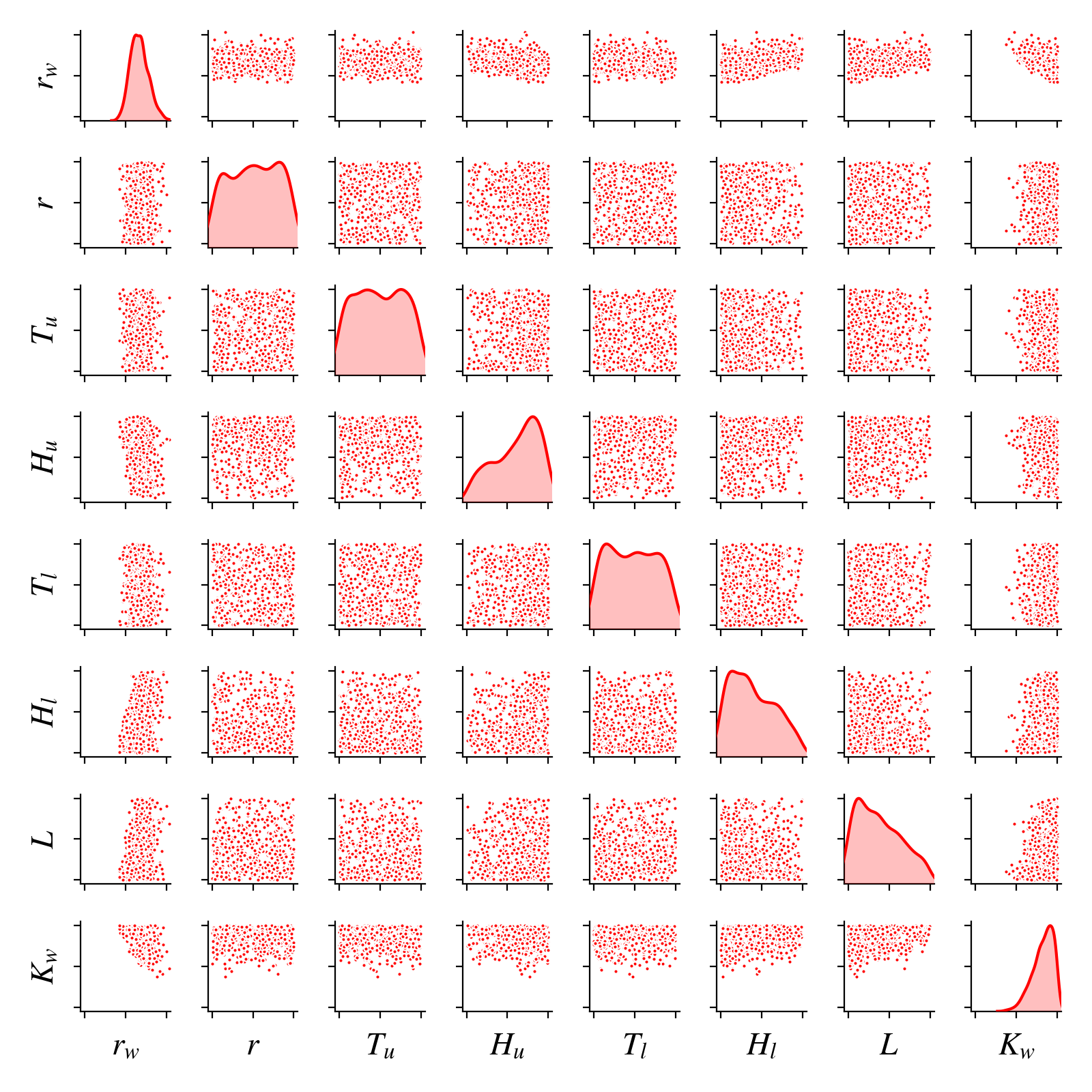}
\caption{Pairwise scatter plot of isolated samples leading to high outputs for the borehole model.}
\label{fig:borehole_top_pdf}
\end{figure}

The results over 30 trials are shown in \cref{fig:borehole_extr}. It can be seen that near high outputs values, $r_w$ has significantly increased importance compared to that predicted by full space Sobol' indices, along with somewhat increased importance for $H_u$, $H_l$ and $L$. For low outputs, the importance ranking is similar to that evaluated from full space Sobol' indices. In both cases, it can be seen that no significant loss of accuracy is incurred by using either polynomial or ridge approximations in place of function evaluations for MCF. In \cref{fig:borehole_bot_pdf,fig:borehole_top_pdf}, the marginal and pairwise distributions obtained from MCF are plotted, showing some correlation between $r_w$ and $K_w$ which is not present in the original full space input distribution. From the skewness indices, we see the increased importance of $r_w$, $H_u$, $H_l$ and $L$ with positive values. This corroborates with the interpretation of skewness indices mentioned in \Cref{sec:ex_analytical}. The agreement is only qualitative however, with the relative ranking of $r_w$ and $K_w$ not in agreement. Again, it is not clear whether the observations can hold generally.


\subsection{Piston cycle time function} \label{sec:piston_ex}
Consider a model of a piston adapted from \cite{kenett2013modern}, which describes the cycle time of the piston as a function of seven parameters. These parameters are shown in \Cref{piston_inputs}. The inputs are independently and uniformly distributed over their support. As in the previous example, we are interested in obtaining the relative importance of each input variable near output extrema. The cycle time in seconds is
\begin{equation}
\tau(\vx) = 2\pi \sqrt{\frac{M}{k+S^2\frac{P_0 V T_a}{T_0 V^2}}},
\end{equation}
where 
\begin{equation}
V = \frac{S}{2k}\left( \sqrt{A^2 + 4k\frac{P_0 V}{T_0} T_a} - A \right)\,
\end{equation}
and
\begin{equation}
A = P_0 S + 19.62M - \frac{kV}{S}.
\end{equation}

\begin{table} [h]
\centering
\begin{tabular}{@{}ccc@{}}
\hline
\textbf{Input} & \textbf{Range} & \textbf{Description} \\
\midrule
$M$ & [30, 60] & Piston mass (kg) \\

$S$ & [0.005, 0.020] & Piston surface area (m$^2$) \\

$V$ & [0.002, 0.010] & Initial gas volume (m$^3$) \\

$k$ & [1000, 5000] & Spring coefficient (N/m) \\

$P_0$ & [90000, 110000] & Atmospheric pressure (N/m$^2$) \\

$T_a$ & [290, 296] & Ambient temperature (K)\\

$T_0$ & [340, 360] & Filling gas temperature (K)\\
\bottomrule
\end{tabular}
\caption{Input variables to the piston model.}
\label{piston_inputs}
\end{table}

\begin{figure}
\centering
\includegraphics[width=1\textwidth]{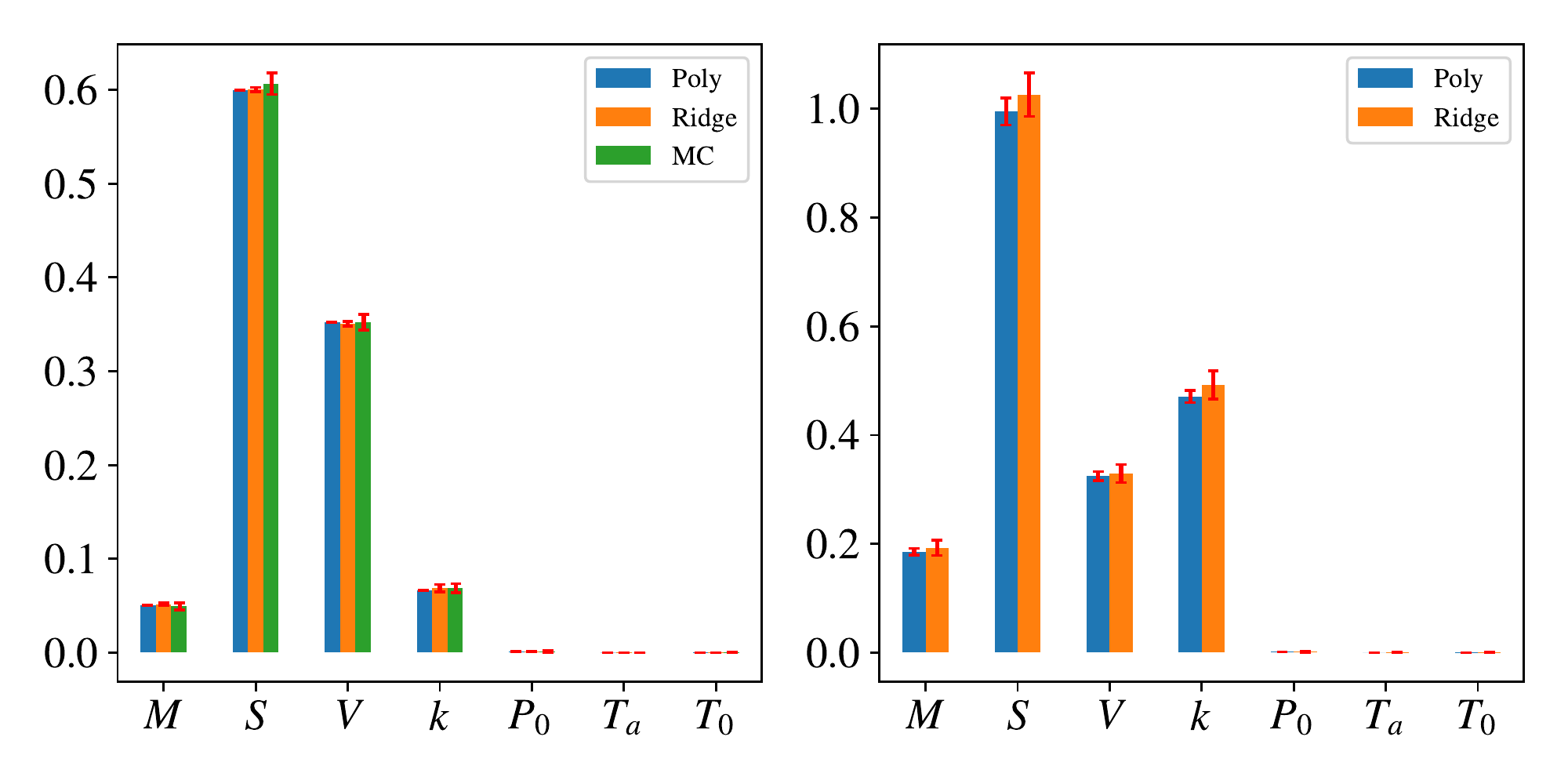}
\caption{Total Sobol' (left) and skewness (right) indices for the piston function.}
\label{fig:piston_inds}
\end{figure}

The total Sobol' and skewness indices are evaluated using the same three methods as in the previous example---Poly, Ridge and MC. For Poly, a degree 4 polynomial on a total order isotropic basis (cardinality 330) is fitted with 500 samples; for Ridge, a degree 4 polynomial over a 4-dimensional subspace is fitted with 300 samples and its coefficients. The results are shown in \cref{fig:piston_inds}, again showing that using ridge polynomial approximations incurs insignificant losses in accuracy. 

\begin{figure}
\centering
\includegraphics[width=1\textwidth]{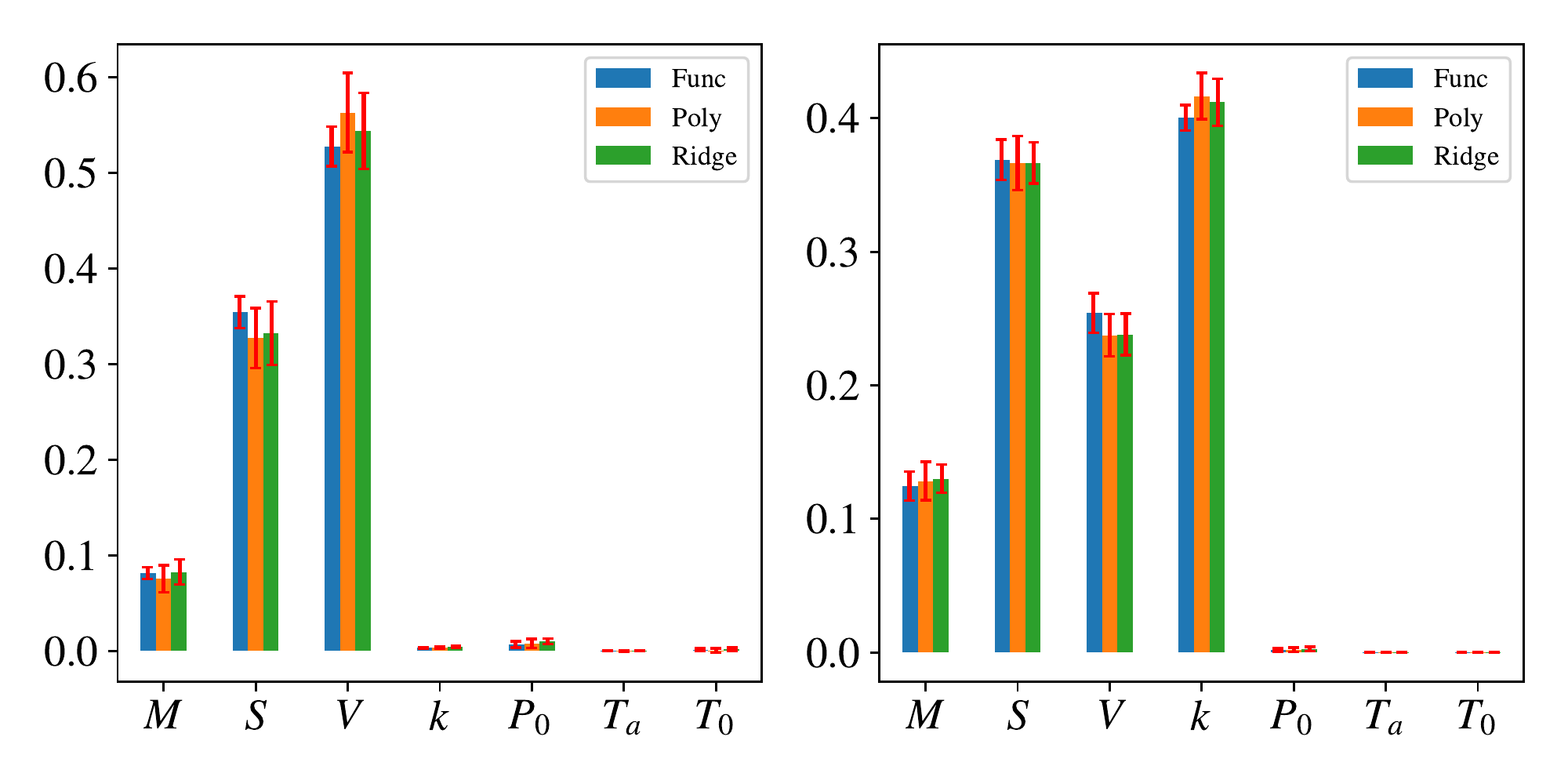}
\caption{Bottom (left) and top (right) Sobol' indices for the piston function using function evaluations (Func), a polynomial approximation (Poly) and a ridge approximation (Ridge).}
\label{fig:piston_extr}
\end{figure}

\begin{figure}
\centering
\includegraphics[width=0.5\textwidth]{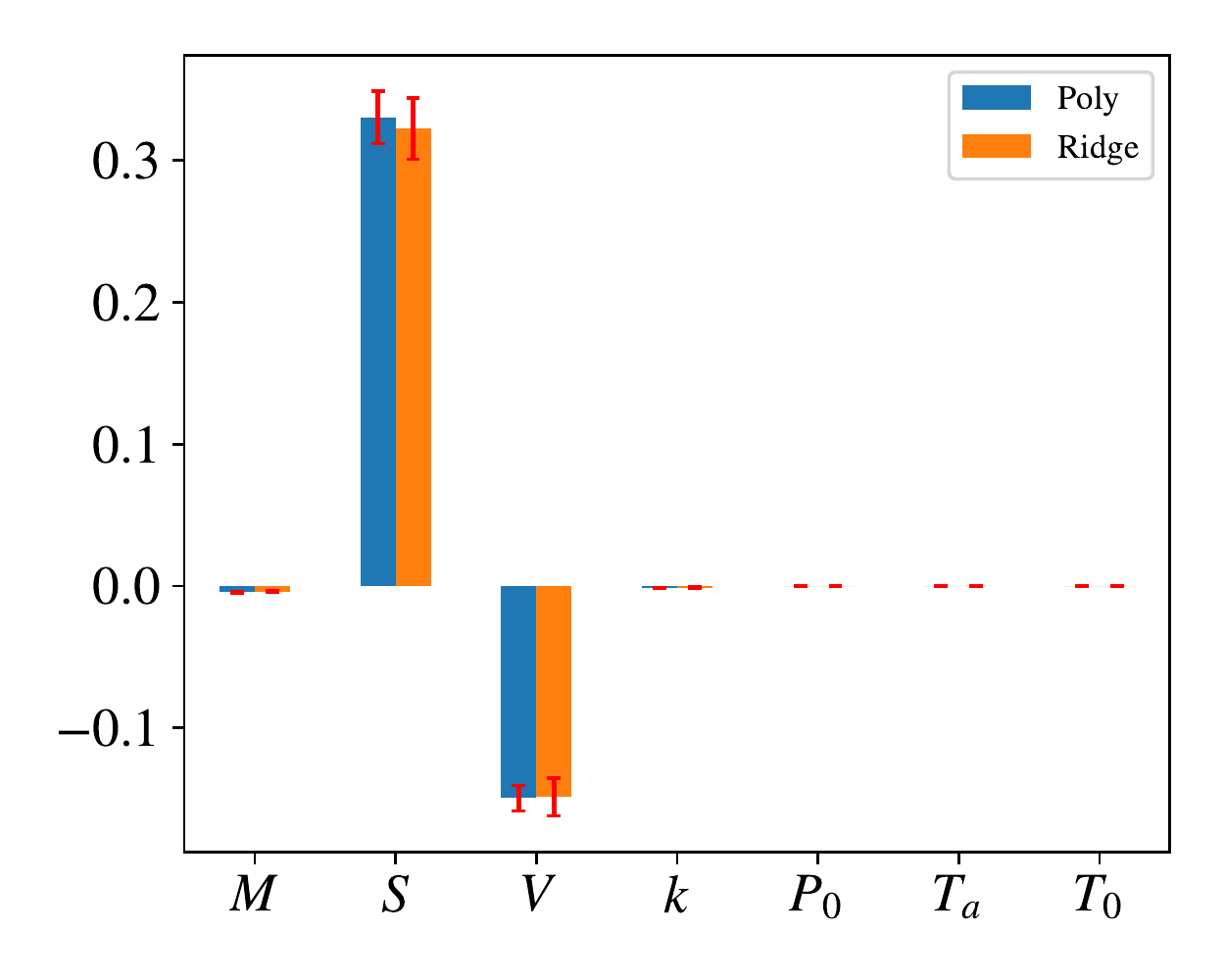}
\caption{First order skewness indices for the piston function.}
\label{fig:piston_first_skew}
\end{figure}

\begin{figure}
\centering
\includegraphics[width=1\textwidth]{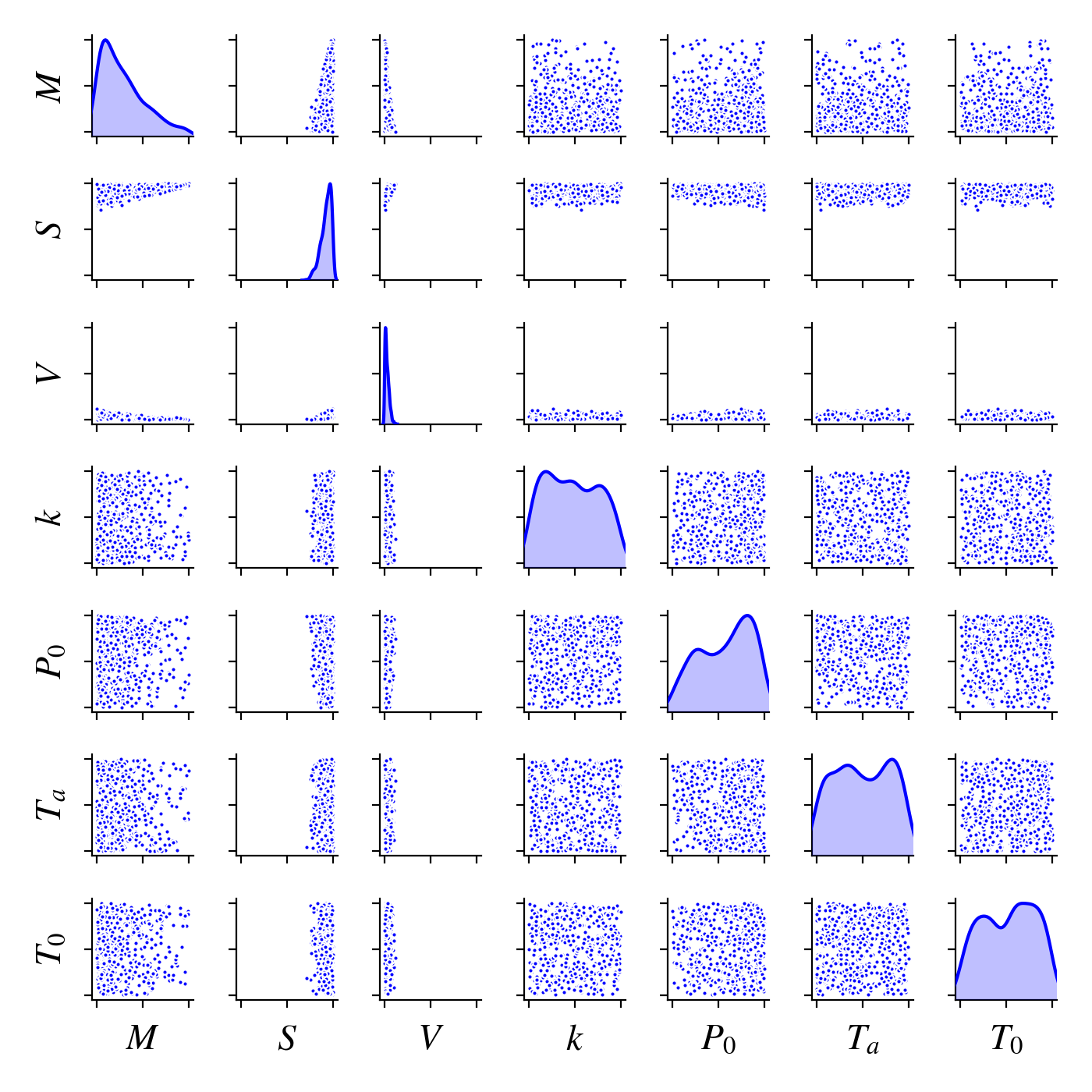}
\caption{Pairwise scatter plot of isolated samples leading to low outputs for the piston model.}
\label{fig:piston_bot_pdf}
\end{figure}

\begin{figure}
\includegraphics[width=1\textwidth]{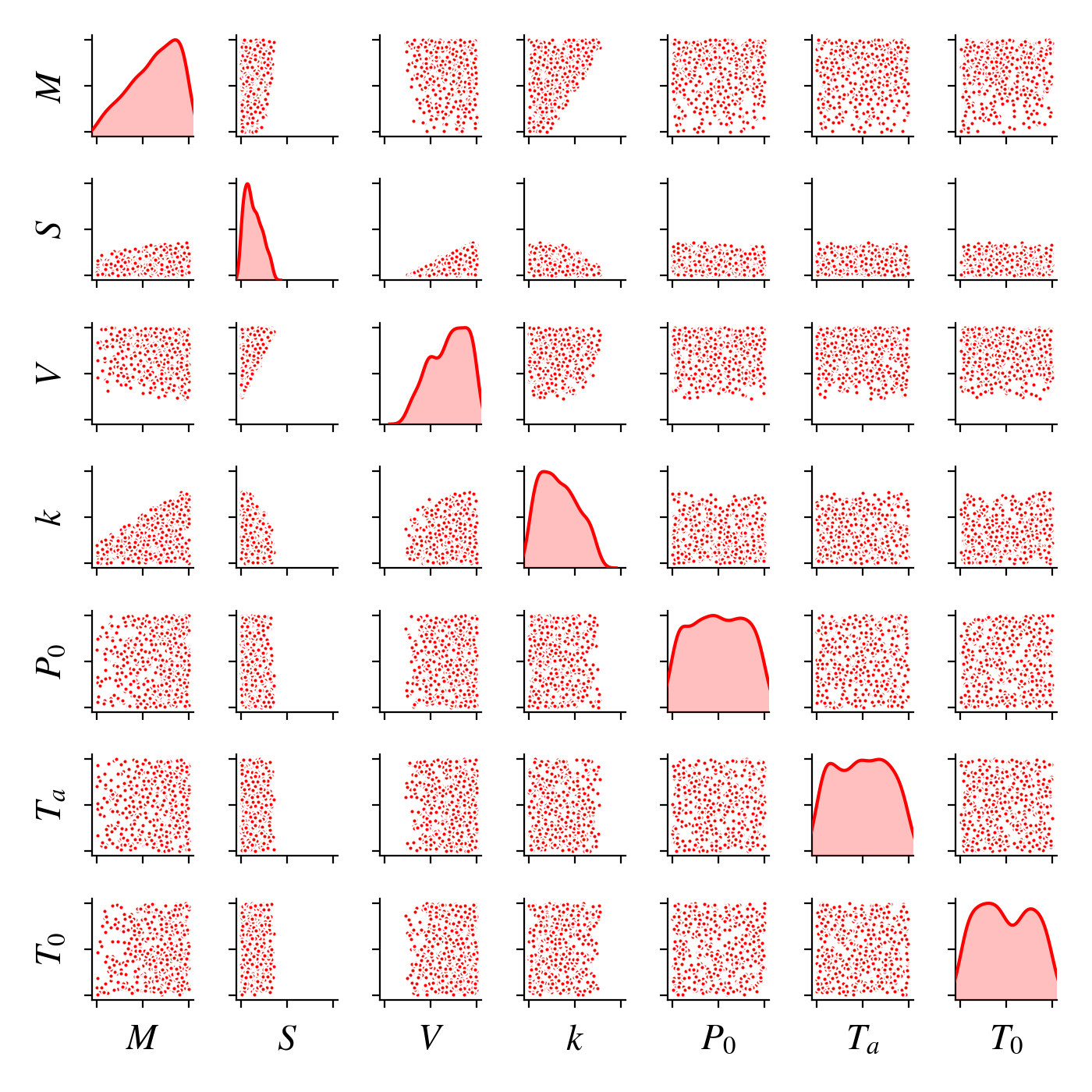}
\caption{Pairwise scatter plot of isolated samples leading to high outputs for the piston model.}
\label{fig:piston_top_pdf}
\end{figure}

The extremum Sobol' indices are evaluated using the same approaches as in the previous example. For Poly, we use a degree 5 polynomial (cardinality 792) trained from 1000 samples; for Ridge, a degree 5 polynomial over a 4-dimensional subspace trained with 700 samples is used. A degree 4 polynomial in the extrema is fitted with an additional 500 training points for evaluating the indices. The results are shown in \cref{fig:piston_extr}. Comparing with full space Sobol' indices, the bottom Sobol' indices show an increased importance of $V$ in driving the output at low values. For high outputs, $k$ becomes important. Examining the total skewness indices, the increased importance at maximum can be seen for $k$, which has a large positive total skewness index.  The information for $V$ is not found in the total skewness indices, but examining the \emph{first order} skewness indices (\Cref{fig:piston_first_skew}), the importance of $V$ at minimum is highlighted by the large negative value. From the top and bottom Sobol' indices, we find that this variable is important throughout the domain including the output extrema at both ends, but it is especially important at the minimum. With skewness indices, it is necessary to examine different orders to draw this conclusion. Again, we note that in the above, empirical observations about skewness indices are described, but no attempt has been made for a unifying theoretical statement---an interesting topic that we invite further discussions about. \Cref{fig:piston_bot_pdf,fig:piston_top_pdf} show the pairwise plots of the extrema PDFs, showing correlations between $S$, $V$ and $k$ near extrema, and justifying the use of polynomials taking this into account.


\section{Conclusions}
In this paper we explore two related threads in sensitivity analysis. Through a polynomial-based approximation framework, we show the possibility of leveraging the ridge structure of a quantity of interest to mitigate the curse of dimensionality and reduce the number of observations required for sensitivity analysis. In addition, the use of polynomials and their ridges is demonstrated in the evaluation of a new type of sensitivity metric that quantifies sensitivities near output extrema---the extremum Sobol' indices, reducing the cost of otherwise expensive Monte Carlo filtering operations. Since some empirical connections can be observed between extremum Sobol' indices and skewness indices, we hope that this work will encourage discussion on the interpretation and utility of skewness indices in extremum sensitivity analysis.

\section*{Acknowledgements}
This work was supported by the Cambridge Trust; Jesus College, Cambridge; the Alan Turing Institute and Rolls-Royce plc. We would also like to thank Art Owen for helpful discussions, and anonymous reviewers for their comments, which helped improve the manuscript.

\appendix
\section{Computing extremum Sobol' indices} \label{sec:extr_sobol}
The algorithm for computing extremum Sobol' indices after obtaining extremum samples from MCF is shown in \Cref{alg:extr_sobols}. Note that Sobol' indices beyond the first order need to be computed recursively by evaluating the associated quantities at lower orders (in line 15).

\begin{algorithm}
\caption{Computing extremum Sobol' indices}
\label{alg:extr_sobols}
\begin{algorithmic}[1]
\Input
\State Extremum distribution $\rho_e(\vx)$, indices of the variables of interest $S$.
\EndInput
\Output
\State Extremum Sobol' index $\sigma_{S}$.
\EndOutput
\State Obtain $N_u$ samples from $\rho_e(\vx)$ and form the Vandermonde matrix
\begin{equation}
\mA_u =  \begin{bmatrix}
\Psi_1 (\mathbf{x}_1) & \Psi_2 (\mathbf{x}_1) & \dots & \Psi_r (\mathbf{x}_1) \\
\Psi_1 (\mathbf{x}_2) & \Psi_2 (\mathbf{x}_2) & \dots & \Psi_r (\mathbf{x}_2) \\
\vdots & \vdots & \ddots & \vdots \\
\Psi_1 (\mathbf{x}_{N}) & \Psi_2 (\mathbf{x}_{N}) & \dots & \Psi_r (\mathbf{x}_{N})
\end{bmatrix}.
\end{equation}
where $(\Psi_i)$ is an orthogonal polynomial basis with respect to the \emph{product of the marginals} of $\rho_e(\vx)$.

\State Perform a QR-decomposition of $(1/\sqrt{N}) \mA_u = \mQ_u \mR$, and form a new polynomial basis $(\Phi_i)$ where
\begin{equation}
\Phi_i(\vx) = \sum_{j=1}^r \Psi_j(\vx) (\mR^{-1})(j, i)
\end{equation}
This basis is now orthogonal with respect to $\rho_e(\vx)$ \cite{jakeman2019polynomial}.
\State Sample $N_e$ points and evaluate the function on these points ($\vf_e$). Solve for the least squares coefficients
\begin{equation}
\underset{\boldsymbol{\alpha}_e }{\text{minimize}} \; \; \left\Vert \mA_{c} \boldsymbol{\alpha}_e - \vf_e \right\Vert _{2},
\end{equation}
where $\mA_{c}$ is the Vandermonde matrix based on $(\Phi_i)$.
\State Initialize empty matrix $B \in \mathbb{R}^{N_u \times r}$.
\For{$i=1,...,r$}
\State Let $\vi$ be the corresponding multi-index of $(\Psi_i)$.
\State $\vi_m = \vi$ with its $s$-th element set to 0 for all $s \in S$.
\State $\vi_v = $ zero vector with its $s$-th element set to that of $\vi$ for all $s \in S$.
\State Set the $i$-th column of $\mB$ as 
\begin{equation}
\mB(:, i) = \mu_{i_m} \mA_u(:, i_v)
\end{equation}
where $\mu_{i_m}$ is the mean of $\mA_u(:, i_m)$. 
\EndFor
\State Calculate $\vm_S = \mB \mR^{-1} \boldsymbol{\alpha}_e - \sum_{T\subset S} \vm_T.$
\State Calculate the empirical covariance $C_S = \text{Cov}[\vm_S, \mA_{c} \boldsymbol{\alpha}_e ]$
\State $\sigma_S = C_S / \sigma$ where $\sigma$ is the output variance, calculated by summing up the squares of the coefficients $\boldsymbol{\alpha}_e$ without the first.
\end{algorithmic}
\end{algorithm}

\section{Computing skewness indices} \label{sec:appen}

As described with \eqref{eqn:skewness_indices} in \Cref{sec:extremeum}, the method for calculating skewness indices via polynomial expansion coefficients is more involved than calculating the Sobol' indices, because the integrals of third-order cross terms cannot be simplified easily with orthogonality. Hence, the expectation terms within each summand need to be evaluated using numerical quadrature. That is, we can approximate the expectation of a function $f(\vx)$ using Gauss quadrature
\begin{equation}
\myEx[f(\vx)] = \int_\mathcal{D} f(\vx) \rho(\vx) d\vx \approx \sum_{q=1}^M f(\vs_q) \nu_q,
\end{equation}
where $(\vs_q,\nu_q)_{q=1}^M$ are the quadrature points and weights. In \Cref{alg:cond_skew} we describe the algorithm implemented in Effective Quadratures \cite{seshadri2017effective-quadratures}. In this algorithm, we have defined the \emph{normalized multi-index} corresponding to the set of variables indexed by $S$ 
\begin{equation}
\bar{\vj} = (\bar{j_1},...,\bar{j_d})
\end{equation}
where $\bar{j_i}=1$ if $i\in S$ and 0 otherwise, and $\odot$ denotes elementwise multiplication. In general, an $O(r^3)$ loop is required for calculating a skewness index, and since $r$ can scale exponentially with the input dimension, this presents a significant computational burden. Our algorithm incorporates several methods to reduce the computational time, namely
\begin{itemize}
\item Using a selection function similar to that proposed in \cite{geraci2016high-order} to avoid computation of terms which we know are zero;
\item Scanning the index set in an $O(r)$ loop before computation to filter out indices that involve more variables than specified. This cuts the cost of computation significantly for indices involving a small number of variables.
\item Using a formulation amenable to matrix operations to avoid explicit loops.
\end{itemize}
The code for this algorithm can be found at \url{https://www.effective-quadratures.org}.

\begin{algorithm}
\caption{Computing conditional skewness terms}
\label{alg:cond_skew}
\begin{algorithmic}[1]
\Input
\State Polynomial basis and coefficients, quadrature points and weights $(\vs_q,\nu_q)_{q=1}^M$, set of variable indices $S$ for desired skewness index.
\EndInput
\Output
\State Skewness index $t_{S}$ with corresponding normalized multi-index $\bar{\vj}$.
\EndOutput
\State Compute the weighted evaluation matrix $\mE_w \in \mathbb{R}^{r\times M}$ such that 
\begin{equation} \label{eqn:ew}
\mE_w = \begin{bmatrix}
\alpha_1 \Psi_1(\vs_1) & \alpha_1\Psi_1(\vs_2) & \dots & \alpha_1\Psi_1(\vs_M)\\
\alpha_2\Psi_2(\vs_1) & \alpha_2\Psi_2(\vs_2) & \dots & \alpha_2\Psi_2(\vs_M)\\
\vdots & \vdots & \vdots & \vdots \\
\alpha_r\Psi_r(\vs_1) & \alpha_r\Psi_r(\vs_2) & \dots & \alpha_r\Psi_r(\vs_M)\\
\end{bmatrix}.
\end{equation}
\State Sum every column to get the polynomial approximant at every point $\tilde{f}(\vs_q)$:
\begin{equation}
\tilde{f}(\vs_q) = \sum_{i=1}^r \alpha_i \Psi_i(\vs_q)
\end{equation}
\State Compute the global skewness $\gamma = \sum_{q=1}^M (\tilde{f}(\vs_q))^3 \nu_q$ via a dot product.
\State Compute the global variance $\sigma = \sum_{\text{all }\vi} \alpha_{\vi}^2$.
\State Initialize set of valid indices $V = \varnothing$.\Comment{Prune indices}
\For{$a=1,...,P$}

\State Compute $ord = \sum_d \bar{a}_d$
		\If{$ord \leq \sum_d \bar{j}_d$}
			\State $V \leftarrow V \cup \{a\}.$
		\EndIf
\EndFor
	\State Set $r_1 = 0$.\Comment{Calculate first sum term}
	\For{$a\in V$}
		\State Compute $s_a = \sum_{q=1}^M (\ve_a^3)_q \nu_q$ via a dot product.
		\State $r_1 \leftarrow r_1 + s_a/(\sigma^{3}\gamma)$.
	\EndFor
	\State Set $r_2 = 0$.\Comment{Calculate second sum term}	
	\For{$a,b\in V$}
		\If {($a_d = 0$ and $b_d \neq 0$ for any $d$) or $a=b$}
			\State Continue.
		\EndIf
		\State Compute $s_{ab} = \sum_{q=1}^M (\ve_a^2 \odot \ve_b)_q \nu_q$ via a dot product.
		\State $r_2 \leftarrow r_2 + 3s_{ab}/(\sigma^{3}\gamma)$.
	\EndFor
\algstore{myalg}
\end{algorithmic}
\end{algorithm}

\begin{algorithm}
\begin{algorithmic}[1]
\algrestore{myalg}
	\State Set $r_3 = 0$.\Comment{Calculate third sum term}	
	\For{$a\in V$, $b\in \{a+1,...,P\} \cap V$, $c=\{b+1,...,P\} \cap V$}
		\If {$\bar{a}_d + \bar{b}_d + \bar{c}_d = 1$ for any $d$}	
			\State Continue.
		\ElsIf{$\bar{a}_d + \bar{b}_d + \bar{c}_d = 2$ and $\bar{a}_d, \bar{b}_d,\bar{c}_d$ are pairwise distinct for any $d$}
			\State Continue.
		\EndIf
		\State Compute $s_{abc} = \sum_{q=1}^M (\ve_a\odot \ve_b\odot \ve_c)_q \nu_q$ via a dot product.
		\State $r_3 \leftarrow r_3 + 6s_{abc}/(\sigma^{3}\gamma)$.
	\EndFor
\State Compute $t_{S} = r_1 + r_2 + r_3$.
\end{algorithmic}
\end{algorithm}

\section{Adaptive least angle regression} \label{sec:appenB}

In this section, the method of adaptive least angle regression (LARS)   implemented for comparison in the numerical examples (\Cref{sec:num_ex}) is briefly described. This method is based on \cite{blatman2011adaptive}, and the interested reader is referred to this work for further details. Define the predictor vector $\phi_i$ to be the evaluation of the $i$-th basis function on the training data set, i.e. the $i$-th column of the polynomial design matrix \eqref{eqn:design_matrix}. The LARS algorithm \cite{efron2004least} proceeds as follows:
\begin{enumerate}
\item Standardize all predictors and start with the residual $\delta = \vy - \bar{\vy}$. Initialize all coefficients $\bal = \mathbf{0}$.
\item Find the predictor $\phi_i$ most correlated with $\delta$.
\item Increase $\alpha_i$ until some other $\phi_j$ has as much correlation with $\vr$ as $\phi_i$ does. 
\item Increase both $\alpha_i$ and $\alpha_j$ towards their joint least squares coefficient on $\{\phi_i,~\phi_j\}$ until another predictor $\phi_k$ has as much correlation with the current residual.
\item Repeat until min($r,~N$) predictors have entered, where $r$ is the number of basis functions and $N$ is the number of available data points.
\end{enumerate}
The adaptive LARS procedure applies LARS for selecting a sparse basis, as follows:
\begin{enumerate}
\item Choose an initial candidate basis with maximum polynomial degree of $p = 1$.
\item Apply LARS to the candidate basis and retain the terms with non-zero coefficients.
\item Calculate an estimate of the error of the model on this pruned basis via the \emph{leave-one-out} estimate on the training data (see \cite{blatman2011adaptive}).
\item Repeat for increasing $p$ and find the model with the least error at order $p^*$.
\item Using the basis found for $p^*$, fit a polynomial using least squares regression.
\end{enumerate}
Note that this procedure is a special case of the method described in \cite{blatman2011adaptive} in two respects: 1) A fixed data set is used instead of a sequential experimental design for fair comparison with other methods in this paper; and 2) An isotropic total order basis is used in all instances, consistent with other methods. For the results in this paper, the version of LARS implemented in the open-source Python package \texttt{scikit-learn} \cite{pedregosa2011scikit-learn} is used.

\section{Additional verification results for the borehole and piston functions} \label{sec:ver}

\begin{figure}
\centering
\includegraphics[width=1\textwidth]{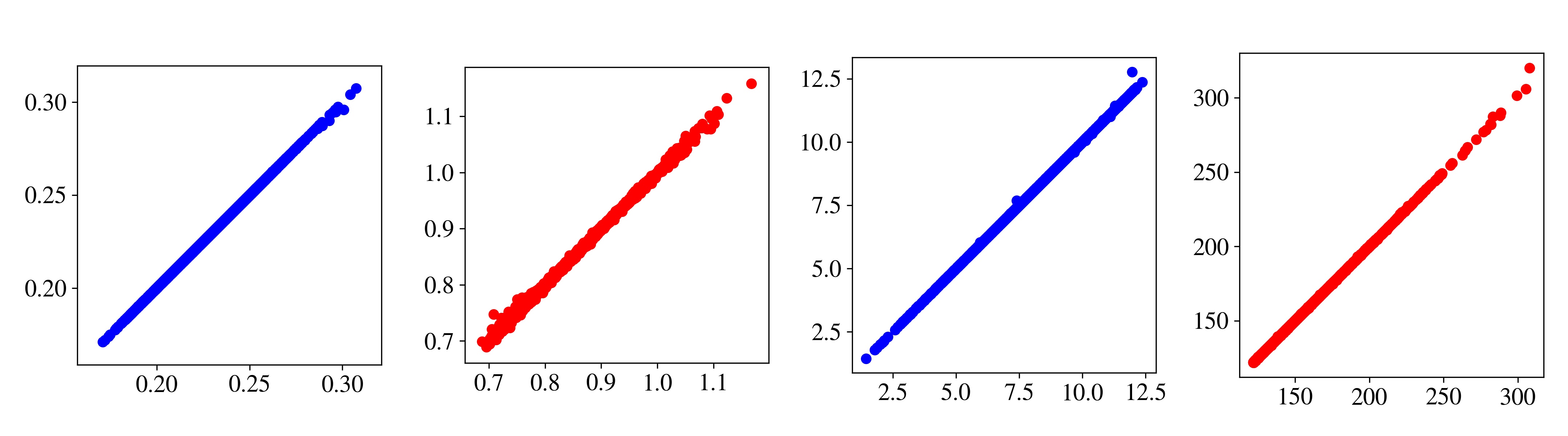}
\caption{Scatter plots of extreme polynomial evaluations against model evaluations for the piston model (left two plots) and borehole function (right two plots). Blue points are for bottom outputs and red points for top outputs.}
\label{fig:poly_verify}
\end{figure}

In this section, the choices of polynomial degrees used for sensitivity analysis of the piston and borehole models are verified. First, \Cref{fig:poly_verify} verifies the accuracy of the extreme polynomials (polynomials fitted to extreme points from MCF) by plotting polynomial evaluations against true model evaluations for points near both extrema for both functions. The points shown are testing points that are not used in the fitting process. This figure shows that the polynomial degree chosen gives sufficient accuracy at the desired output range. 

\begin{figure}
\centering
\includegraphics[width=1\textwidth]{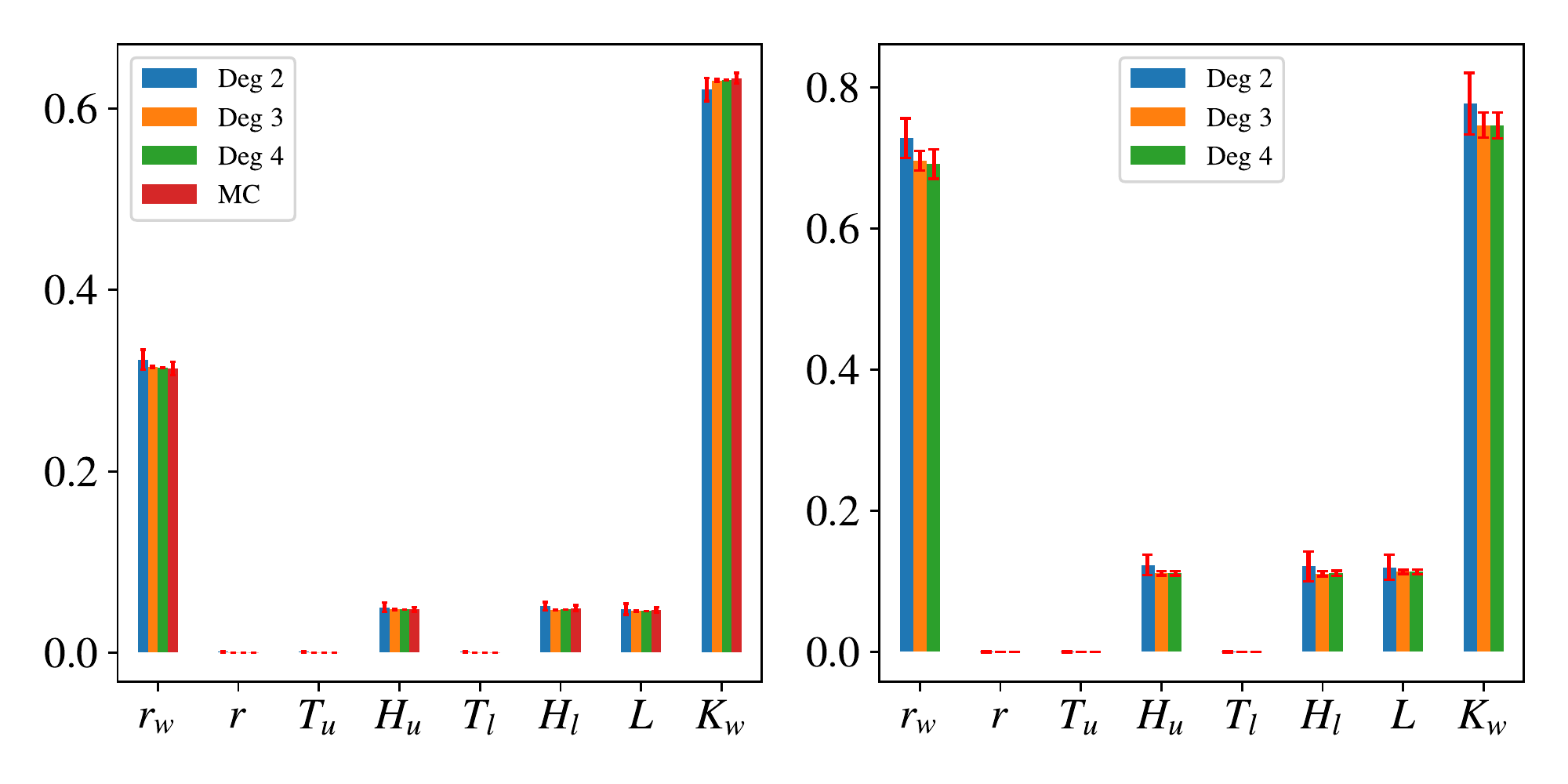}
\caption{Total Sobol' and skewness indices for the borehole model with polynomial approximations at various degrees.}
\label{fig:borehole_inds_verify}
\end{figure}

\begin{figure}
\centering
\includegraphics[width=1\textwidth]{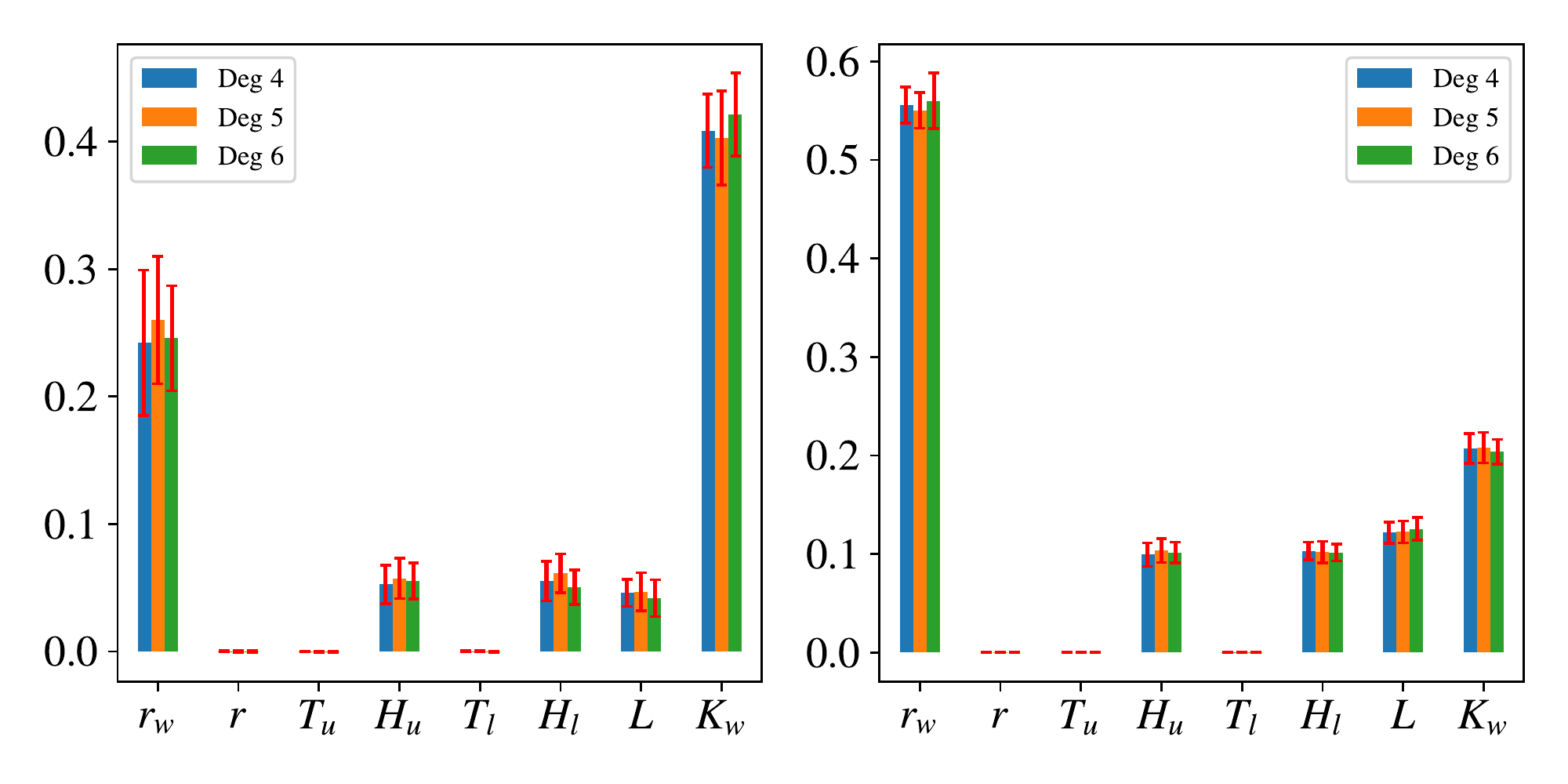}
\caption{Extremum Sobol' indices for the borehole model with polynomial approximations at various degrees.}
\label{fig:borehole_extr_verify}
\end{figure}

\begin{figure}
\centering
\includegraphics[width=1\textwidth]{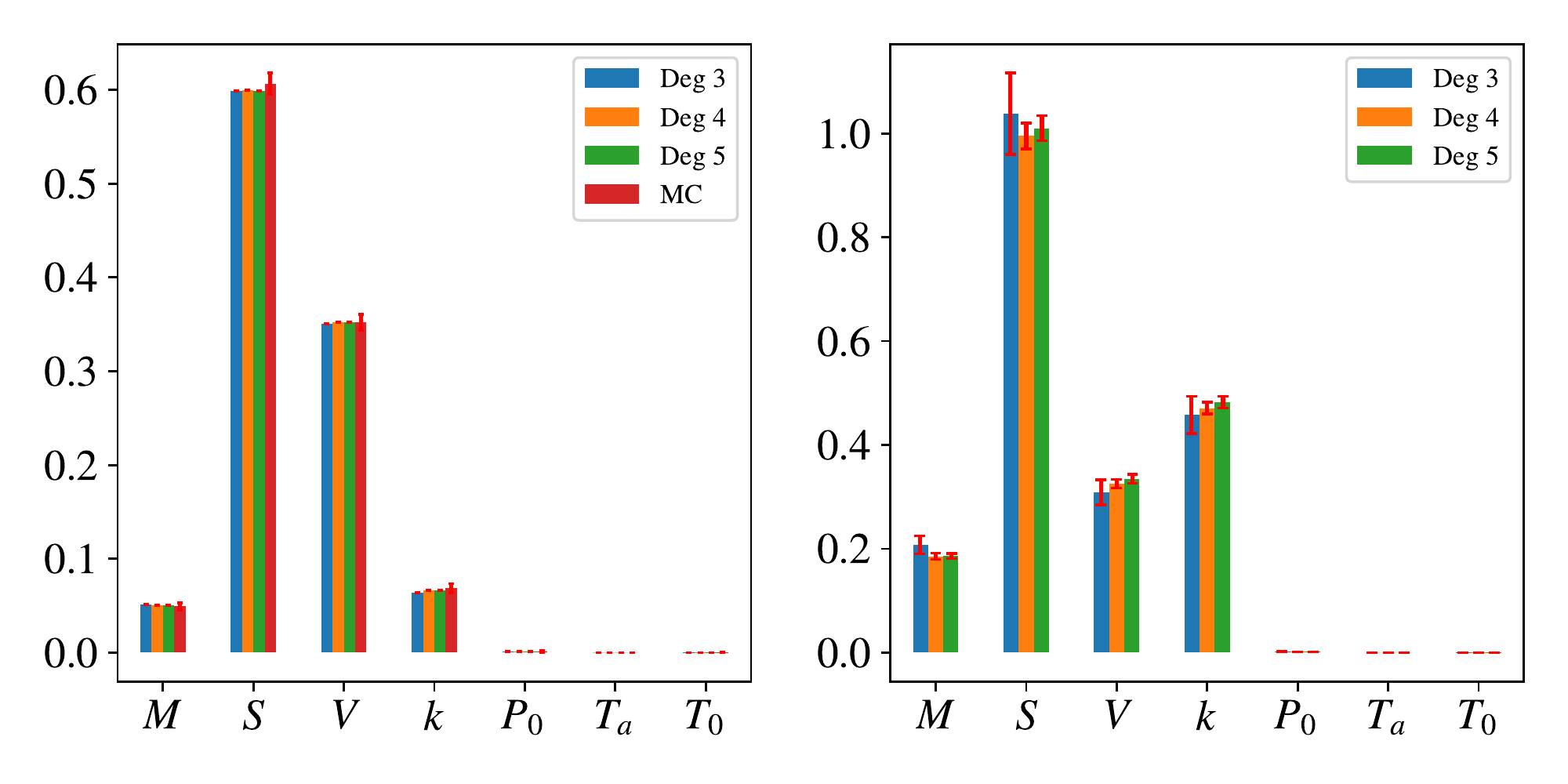}
\caption{Total Sobol' and skewness indices for the piston model with polynomial approximations at various degrees.}
\label{fig:piston_inds_verify}
\end{figure}

\begin{figure}
\centering
\includegraphics[width=1\textwidth]{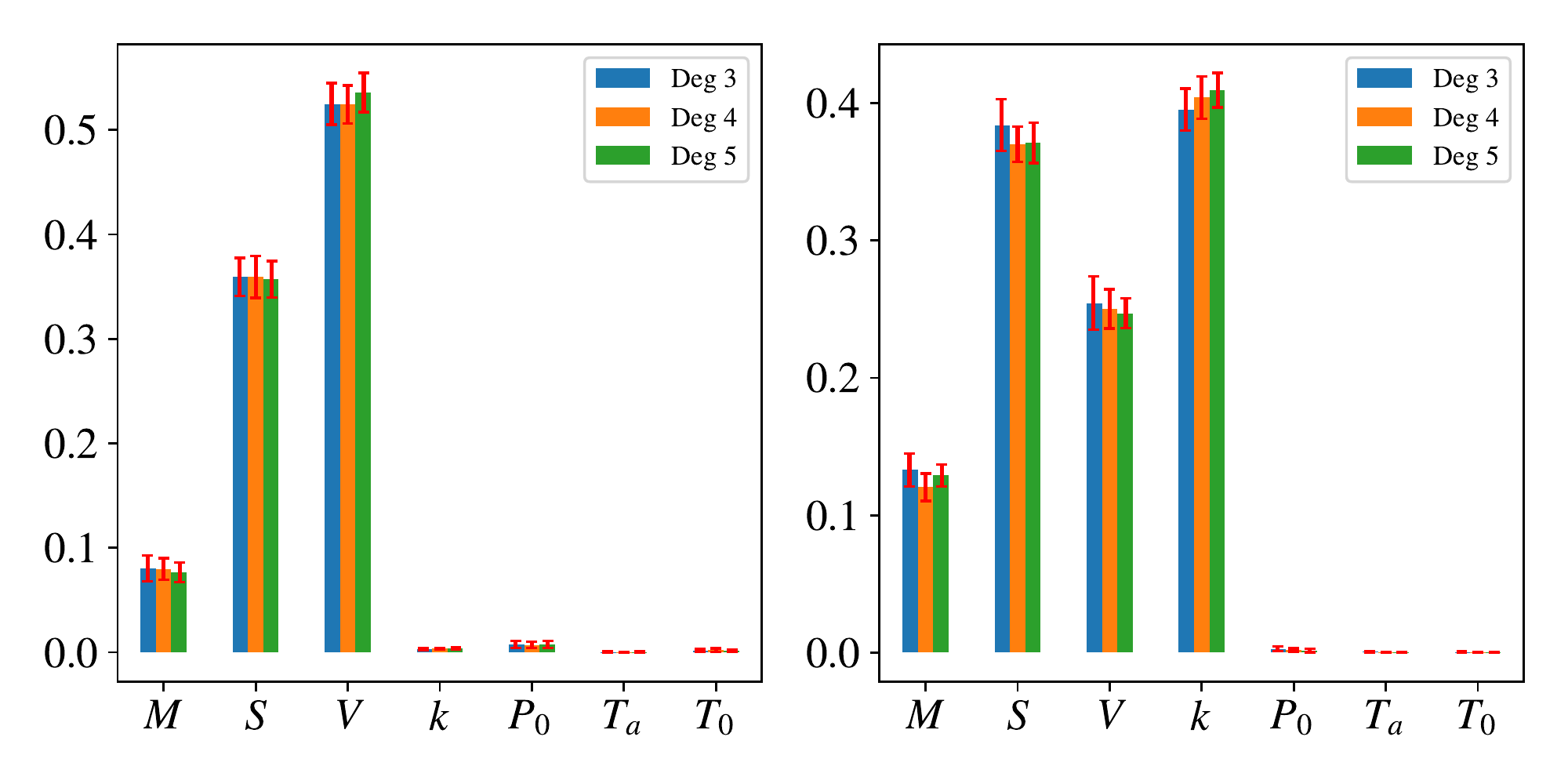}
\caption{Extremum Sobol' indices for the piston model with polynomial approximations at various degrees.}
\label{fig:piston_extr_verify}
\end{figure}

Second, the sensitivity indices for the borehole and piston functions are computed again with perturbations to the degrees of the polynomial approximations. Note that for extremum Sobol' indices, the polynomial degree refers to that of the polynomial fitted to the extremum samples (according to the extremum distribution), not the polynomial used to obtain MCF samples. The latter is already compared with actual function evaluations which ensures its accuracy. \Cref{fig:borehole_inds_verify,fig:borehole_extr_verify,fig:piston_inds_verify,fig:piston_extr_verify} show that perturbing the polynomial degrees produces negligible changes to the values of the indices, thus confirming that the degrees used in the main text are appropriate.

\bibliography{ress}

\end{document}